\documentclass[preprint,12pt]{elsarticle}

\usepackage{subcaption}
\usepackage{graphicx}
\usepackage{booktabs}
\usepackage{geometry}
\usepackage{color}
\usepackage{stmaryrd}
\usepackage{algorithm}
\usepackage{algpseudocode}
\usepackage{float}

\usepackage{amssymb}
\usepackage{amsmath}
\usepackage{amsthm}

\theoremstyle{definition}
\newtheorem{example}{Example}[section]

\journal{Journal of Scientific Computing}

\begin{document}

\begin{frontmatter}



\title{PINN-Based Kolmogorov-Arnold Networks with RAR-D Adaptive Sampling for Solving Elliptic Interface Problems} %% Article title


\author[label1]{Zijuan Xin} %% Author name
\author[label1]{Chenyao Wang}
\author[label1]{Feng Shi}
\author[label2]{Yizhong Sun \corref{cor1}}

\affiliation[label1]{organization={College of Science, Harbin Institute of Technology},%Department and Organization
            %addressline={}, 
            city={Shenzhen},
            %postcode={}, 
            %state={},
            country={China}}

\affiliation[label2]{organization={Department of Mathematics, Hong Kong Baptist University},%Department and Organization
            %addressline={}, 
            city={Hong Kong},
            %postcode={}, 
            %state={},
            country={China}}

\cortext[cor1]{Corresponding author (email: bill950204@126.com). }

%% Abstract
\begin{abstract}
Physics-Informed Neural Networks (PINNs) have become a popular and powerful framework for solving partial differential equations (PDEs), leveraging neural networks to approximate solutions while embedding PDE constraints, boundary conditions, and interface jump conditions directly into the loss function. However, most existing PINN approaches are based on multilayer perceptrons (MLPs), which may require large network sizes and extensive training to achieve high accuracy, especially for complex interface problems. In this work, we propose a novel PINN architecture based on Kolmogorov–Arnold Networks (KANs), which offer greater flexibility in choosing activation functions and can represent functions with fewer parameters. Specifically, we introduce a dual KANs structure that couples two KANs across subdomains and explicitly enforces interface conditions. To further boost training efficiency and convergence, we integrate the RAR-D adaptive sampling strategy to dynamically refine training points. Numerical experiments on the elliptic interface problems yield more uniform error distributions across the computational domain, which demonstrates that our PINN-based KANs achieve superior accuracy with significantly smaller network sizes and faster convergence compared to standard PINNs. 
\end{abstract}


%\item Introduce a dual-subdomain PINN-Based Kolmogorov--Arnold Network for elliptic interface problems.
%\item Enforce interface jump and boundary conditions explicitly through a coupled physics-informed loss.
%\item Integrate RAR-D adaptive sampling to refine collocation points near interface-dominated residual regions.
%\item Achieve higher accuracy and faster convergence than standard MLP-based PINNs.
%\item Produce more uniform error distributions across the computational domain for representative elliptic interface benchmarks.

%% Keywords

\end{frontmatter}



\section{Introduction}
\label{sec1}

%% Labels are used to cross-reference an item using \ref command.
Elliptic interface problems play a central role in computational science and engineering, with applications in multiphase fluid flow \cite{1}, electromagnetics \cite{6}, and heat conduction across heterogeneous materials \cite{4}. A typical model involves discontinuous coefficients and interface jump conditions, which may cause discontinuities in the solution or its flux across an embedded interface. These non-smooth features make it difficult to obtain high accuracy using standard discretizations unless the interface is carefully handled. Traditional approaches, including finite element \cite{8}, immersed methods \cite{Ji_immersed_JSC} and meshless methods \cite{EABE_Mu2024}, have achieved strong performance. Recent progress on interface-unfitted discretizations and efficient solvers further highlights the continuing interest in robust interface treatments \cite{Chu_MG_IFE_JSC}.  However, for complex interface geometry, they often require mesh construction or special treatments near the interface, which can increase implementation complexity and computational cost and may affect stability.

Recently, deep learning has been used to build meshfree solvers for partial differential equations (PDEs). Physics-Informed Neural Networks (PINNs) \cite{14,Cuomo_PINNs_JSC} approximate the solution by a neural network and enforce the PDE and boundary conditions through a loss function, using automatic differentiation and GPU acceleration for efficient training. This framework is attractive for irregular domains and can reduce the burden of mesh generation. Motivated by domain decomposition, deepDDM \cite{18, JSC_DDMPINNs} trains separate networks in subdomains and couples them through interface constraints for elliptic problems. Related decomposition-based PINN variants include cPINNs \cite{19}, which impose flux conservation across subdomains, XPINNs \cite{21}, which extend decomposition ideas to broader nonlinear PDE settings, and I-PINNs \cite{20}, which explicitly incorporate interface jump conditions for interface problems.

Despite these developments, elliptic interface problems remain challenging for standard PINNs. Most existing PINN-type methods use multilayer perceptrons (MLPs) \cite{22} as the basic network model. Although MLPs have the universal approximation property \cite{24}, elliptic interface solutions are typically non-smooth in the neighborhood of the interface, and sharp derivative/flux jumps are difficult to resolve efficiently. In practice, errors and residuals often become large near interfaces, so achieving high accuracy may require large networks and extensive training. Closely related observations have also been reported in discontinuity computing with PINNs, where non-smooth structures can trigger pronounced training difficulties unless additional mechanisms are introduced \cite{Liu_Discontinuity_JSC}.  In addition, PINN performance depends strongly on the placement of collocation points: uniform sampling may place too few points near the interface where the residual is most significant, which slows convergence and yields non-uniform errors.

Two lines of improvement are therefore natural for interface problems: (i) using a more flexible network backbone to better fit interface-dominated solution features with fewer parameters, and (ii) using adaptive sampling to place more training points where the residual indicates insufficient accuracy. For the network backbone, Kolmogorov–Arnold Networks (KANs) \cite{27}, motivated by the Kolmogorov–Arnold representation theorem \cite{28,29}, provide an alternative to MLPs. Unlike MLPs that rely on fixed activation functions, KANs use learnable one-dimensional functions (often represented by B-splines) on network edges \cite{27}, allowing the nonlinear components to adapt during training. KANs have been explored in operator learning and physics-informed modeling, including KAN-based DeepONets \cite{30} and KAN-related PINN frameworks \cite{31,33,35}. For adaptive sampling, residual-based strategies have been systematically studied and shown to improve robustness for problems with sharp or localized features \cite{40,41,42}. In particular, Residual-based Adaptive Refinement with Diversity (RAR-D) \cite{40} refines high-residual regions while maintaining coverage of the full domain, which is especially relevant when large residuals concentrate near interfaces.

In this work, we propose PINN-based KANs with a dual-network domain decomposition structure and RAR-D adaptive sampling for solving elliptic interface problems. The computational domain is partitioned into subdomains, and each subdomain is represented by a dedicated KAN. The subdomain networks are coupled by explicitly enforcing interface jump conditions in the loss function. We further integrate RAR-D \cite{40} to adaptively update collocation points based on the evolving residual distribution, improving convergence and reducing interface-dominated errors. Numerical experiments demonstrate that the proposed PINN-based KANs achieve more uniform error distributions and reach higher accuracy with smaller network sizes and faster convergence than standard MLP-based PINNs.

The paper is organized as follows. Section 2 introduces the notation and a model elliptic interface problem, and reviews a baseline dual-PINNs formulation with standard MLP backbones for enforcing the PDE, boundary conditions, and interface jump conditions. Section 3 presents our proposed PINN-based KANs framework and its dual-network coupling across subdomains, and then introduces the RAR-D adaptive sampling strategy for improving training efficiency in interface-dominated regimes.  Section 4 presents several numerical examples, to show and compare the accuracy and interpretability of PINN-based KANs and standard PINNs.

\section{Model Problem and Baseline Dual-PINNs for Interface Problems}%Problem Setup}
\subsection{Elliptic interface problems}
Firstly, we present the model problem and the corresponding interface conditions. Assume that $\Omega$ is a bounded domain in $R^{2} $, with the boundary $\partial \Omega$ and the interface $\Gamma$, by which $\Omega$ is decomposed into two non-overlapping subdomains $\Omega _{1}$ and $\Omega _{2}$.  The elliptic interface problem model used in this paper can be described as
\begin{equation}\label{eq:model}
\begin{aligned}
  -\nabla \cdot\left(a_{i} \nabla u_{i}\right) &= f_{i},  && \text{in } \Omega_{i}, ~ i=1, 2,  \\
  \llbracket a \nabla u \cdot \mathbf{n} \rrbracket &= \psi,  && \text{on } \Gamma,  \\
  \llbracket u \rrbracket &= \varphi,  && \text{on } \Gamma,   \\
  u_{i} &= g_{i},  && \text{on } \partial \Omega_{i} \backslash \Gamma,
\end{aligned}
\end{equation}
where $u_i=u_i(x,y)$ is the solution in $\Omega_i$,  $f_i$ represents the source term, and $\mathbf{n}$ denotes the unit normal on $\Gamma$. %Assume that $\partial \Omega$ is a Dirichlet boundary, and the interface $\Gamma$ is a closed curve or surface.
We assume that both $\partial \Omega$ and $\Gamma$ are Lipschitz continuous.
The jump of $u$ across the interface is denoted by $\llbracket u \rrbracket: =u|_{\Omega_2}-u|_{\Omega_1}$. The coefficients
\begin{equation}
a(\mathbf{x}) = \left\{
\begin{aligned}
a_{1}(\mathbf{x}), ~& \text { if }~ \mathbf{x} \in \Omega_{1}, \\
a_{2}(\mathbf{x}), ~& \text { if }~ \mathbf{x} \in \Omega_{2},
\end{aligned}
\right.
\end{equation}
are piecewise spatial functions. 
%The unknown part of this problem is the exact solution $u*$, while others are given in advance.

To solve \eqref{eq:model} in a meshfree manner, we approximate the subdomain solutions by neural networks and enforce the governing equations, boundary conditions, and interface jump conditions through a physics-informed least-squares loss. For clarity and for later comparison, we first summarize a standard dual-MLP PINN formulation for elliptic interface problems, which will serve as a baseline in our numerical experiments.

\subsection{Dual-MLP PINNs for elliptic interface problems}
Deep Neural Networks (DNNs) are models based on multi-layer neuron structures, possessing powerful capabilities of data processing and representation. DNNs usually consist of multiple hidden layers, where each neuron in a given layer is connected to neurons in both the previous and subsequent layers. By enforcing weighted connections and non-linear activation functions, DNNs can model complex mapping relationships. DNNs have demonstrated exceptional performance in various fields, such as image recognition and natural language processing, establishing themselves as core methodologies in modern artificial intelligence research\cite{48}. %We connect each neuron in every layer. 
More specifically, for a neural network with $k$ layers, the output can be written as:
\begin{align}
N^k=\mathbf{W}^k\sigma(N^{k-1}(\mathbf{x}))+\mathbf{b}^k, \quad2\leq k\leq L, 
\end{align}
where $\mathbf{x}$ is the input variable of neural network, $\sigma$ stands for the activation function. The final layer, $N^L(\mathbf{x})=\mathbf{W}^L N^{L-1}(\mathbf{x}) +\mathbf{b}^L$, employs an identity activation function. Let $\theta=\{\mathbf{W}^{k},\mathbf{b}^{k}\}\in\mathrm{V}$ represent the set of all weights and biases, where $\mathrm{V}$ denotes the parameter space. 

However, traditional neural network methods approximate the solution $u(\mathbf{x})$ over the entire domain $\Omega$ using a single DNN, often resulting in a slow optimization process and difficulty in achieving sufficient accuracy. By partitioning $\Omega$  into two disjoint subdomains, $\Omega _{1}$ and $\Omega _{2}$, based on the interface, we employ two separate DNN structures in $\Omega _{1}$ and $\Omega _{2}$ to approximate $u|_{\Omega_1}$ and $u|_{\Omega_2}$ denoted as $u_i(\mathbf{x}), i=1, 2, \mathbf{x}\in\Omega_i$.Thus, we have
\begin{align}
u_i(\mathbf{x})\approx U_{i, N}(\mathbf{x}, \theta)=N_i^L\circ N_i^{L-1}\circ\cdots\circ N_i^2\circ N_i^1(\mathbf{x}), \quad i=1, 2.  
\end{align}

%\subsection{Physics-Informed Neural Networks}

In \cite{14}, the authors propose a deep neural network to approximate the solution of partial differential equations, referred to as the 
u-network. Automatic differentiation is then applied to obtain the differential operators of the equation, resulting in an f-network that captures the equation’s physical information. Subsequently, the boundary function and internal loss function are formulated based on the least squares principle.

After determining the inputs of the neural network, we need to use the given boundary conditions and equation information to construct the loss function. Typically, the least squares method is employed, along with automatic differentiation techniques \cite{45}. When solving the elliptic interface problem, we use the random Latin hypercube method to extract the data points, and need to establish two separate PINNs to approximate $u_{1}$ and $u_{2}$, respectively. The total loss of the elliptic interface problem is written as
\begin{align}
\mathcal{L}(\theta) = L_{\Omega_{1}} + L_{\Omega_{2}} + L_{\Gamma} + L_{\partial\Omega_{1}} + L_{\partial\Omega_{2}},
\end{align}
where $L_{\Omega_{i}} \ (i=1,2)$ are the loss of the residuals for the governing equations with respect to the subdomains $\Omega_{i}$, and $L_{\partial\Omega_i}$ and $L_{\Gamma}$ represent the loss associated with the boundary of corresponding domain $\Omega_i$ and the interface conditions, namely
\begin{align}
L_{\Omega_{i}} &= \frac{1}{N_{i}} \sum_{j=1}^{N_{i}}
\left| -\nabla \cdot (a_{i} \nabla u_{i}(\mathbf{x}^{j};\theta_i)) - f_{i}(\mathbf{x}^j) \right|^{2}, \quad i=1,2,\\
L_{\Gamma} &= \frac{1}{N_{\Gamma}} \sum_{j=1}^{N_{\Gamma}} \left|\llbracket u(\mathbf{x}_{\Gamma}^{j};\theta) \rrbracket - \varphi \right|^{2} +  \frac{1}{N_{\Gamma}} \sum_{j=1}^{N_{\Gamma}} \left|\llbracket a \nabla u(\mathbf{x}_{\Gamma}^{j};\theta) \cdot \mathbf{n} \rrbracket - \psi \right|^{2},   \\
L_{\partial\Omega_{i}} &= \frac{1}{N_{\partial\Omega_i}} \sum_{j=1}^{N_{\partial\Omega_i}}
\left| u_{i}(\mathbf{x}_{\partial\Omega_i}^{j};\theta_i) - g_{i}(\mathbf{x}_{\partial\Omega_i}^{j}) \right|^{2}, \quad i=1,2.
\end{align}
where $\{\mathbf{x}^j\}_{j=1}^{N_{i}}$ denotes randomly sampled points within the subdomains $\Omega_{i}$ ($i=1,2$). $\{\mathbf{x}_{\Gamma}^{j}\}_{j=1}^{N_{\Gamma}}$ and $\{\mathbf{x}_{\partial \Omega_i}^{j}\}_{j=1}^{N_{\partial \Omega_i}}$ represent the interface and boundary points, respectively.
%$N_{1}, N_{2}, N_{\Gamma}, N_{\partial\Omega}$ indicate the total numbers of residual points for the two subdomains, boundary points, and interface points.

\begin{figure}[h]
    \centering
    \includegraphics[width=0.95\textwidth]{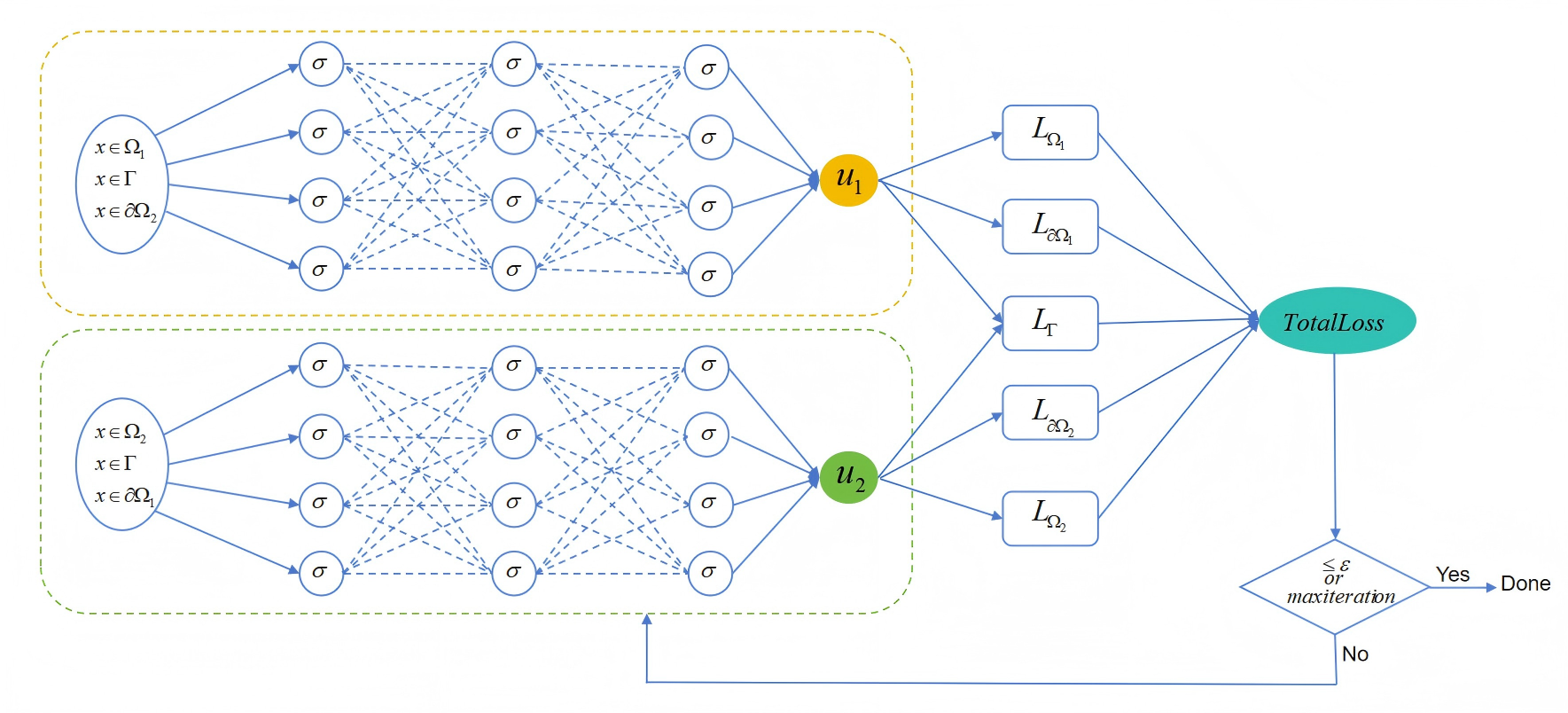} % 替换为图片的实际路径
    \caption{A dual PINNs structure for solving elliptic interface problem.}
    \label{fig:Schematic_PINNS}
\end{figure}

To solve the elliptic interface problem, we employ a dual PINNs architecture to approximate the solution $u_i(\mathbf{x})$. The dual PINNs structure to solve the elliptic interface problem is shown in Figure \ref{fig:Schematic_PINNS}. The basic idea is to use two neural networks to approximate the solution of the elliptic partial differential equation, where each network is responsible for one subdomain of the problem. The networks are coupled through the interface conditions, which enforce the continuity of both the solution and its flux across the interface.

To achieve this, a physics-informed loss function is constructed, which is based on the residuals of the given partial differential equation as well as the interface and boundary conditions, as follows:
\begin{align}
\theta^{*}=\underset{\theta}{\operatorname{argmin}} \mathcal{L}(\theta),
\end{align}
where $\theta^{*}$ is the minimizer of the total loss $\mathcal{L}(\theta)$.

Specifically, the residual terms are derived from substituting the neural network outputs into the PDE, while the interface terms enforce continuity of both the solution and its normal derivative at the interface. The boundary terms enforce the prescribed Dirichlet or Neumann conditions on the boundary of the computational domain.

%Compared to traditional numerical methods, this dual-network approach has the advantage of naturally handling discontinuous coefficients across interfaces, making it well-suited for problems with complex interface geometries. Naturally, this motivates the introduction of the following dual KANs architecture.

The above dual-MLP PINN provides a standard baseline for elliptic interface problems. In our experiments, however, we observe that accurately resolving interface-induced non-smooth features may require large network capacity and careful point placement near the interface. This motivates the PINN-based KANs formulation and the adaptive sampling strategy developed in the next section.

\section{PINN-Based KANs with RAR-D Adaptive Sampling }

Building on the baseline dual-PINN formulation in Section 2, we now present our proposed method. The key idea is to replace the MLP backbone by Kolmogorov--Arnold Networks (KANs) in each subdomain and to couple the subdomain networks through explicit interface losses. We then introduce an adaptive collocation strategy based on RAR-D to improve training efficiency when the residual concentrates near the interface.

\subsection{Kolmogorov-Arnold Networks}

The Kolmogorov--Arnold representation theorem states that a continuous multivariate function can be represented by superpositions of continuous univariate functions \cite{28,29}. Motivated by this theorem, Kolmogorov--Arnold Networks (KANs) parameterize each edge by a learnable univariate function, typically implemented via B-splines \cite{27}. This design enhances the expressivity of the network, enabling it to approximate complex multivariate functions with fewer layers or parameters compared to conventional architectures.

Specifically,the Kolmogorov–Arnold theorem ensures that for any continuous function $f:[0,1]^n\to\mathbb{R}$, there exist continuous univariate functions $\Phi_{q}$ and $\phi_{q, p}$ such that
\begin{align}
f\left(x_{1}, x_{2}, \ldots, x_{n}\right)=\sum_{q=1}^{2 n+1} \Phi_{q}\left(\sum_{p=1}^{n} \phi_{q, p}\left(x_{p}\right)\right).
\end{align}
KANs implement this decomposition by introducing KANLinear layers, which generalize the standard linear layers of neural networks by learning these univariate functions during training.

A general KAN network is a composition of $L$ layers, namely given an input vector $\mathbf{x}\in\mathbb{R}^{n_0}$,  the output of KAN is
\begin{align}
\mathrm{KAN}(\mathbf{x})=(\boldsymbol{\Phi}^{L}\circ\boldsymbol{\Phi}^{L-1}\circ\cdots\circ\boldsymbol{\Phi}^2\circ\boldsymbol{\Phi}^1)(\mathbf{x}),
\end{align}
where $\mathbf{\Phi}^k$$(2\le k\le L)$ is the function matrix corresponding to the $k$-th KAN layer. Each layer applies transformations based on univariate functions, defined as
\begin{align}
\mathbf{\Phi}^k\left(\mathbf{x}^{(k)}\right) = 
\begin{pmatrix}
\phi_{k,1,1}(\cdot) & \cdots & \phi_{k,1,n_k}(\cdot) \\
\phi_{k,2,1}(\cdot) & \cdots & \phi_{k,2,n_k}(\cdot) \\
\vdots & \ddots & \vdots \\
\phi_{k,n_{k+1},1}(\cdot) & \cdots & \phi_{k,n_{k+1},n_k}(\cdot)
\end{pmatrix}
(\mathbf{x}^{(k)}),
\end{align}
where $n_k$ is the number of input nodes for the $k$-th layer and $\phi_{k,i,j}$ is the $k$-th layer’s univariate activation function, connecting its $i$-th input node to its $j$-th output node in the network’s computational graph. Obviously, the output nodes of one layer serve as the input nodes of the subsequent layer. In the case of Kolmogorov–Arnold Networks (KANs), however, the notion of a “layer” is not defined by the collection of nodes, but rather by the edges connecting them, where the activation functions are located. Specifically, for a layer with $n_k$ input nodes and $n_{k+1}$ output nodes, the  total number of univariate activation functions is given by the product $n_k\cdot n_{k+1}$. The architecture of a KAN can be fully specified by a sequence of integers, where each pair of adjacent entries represents the input and output dimensions of a particular layer. In general, the architecture of an $L$-layered KAN is written as $[n_0,n_1,\ldots,n_L]$.

\begin{figure}[h]
    \centering
    \includegraphics[width=0.95\textwidth]{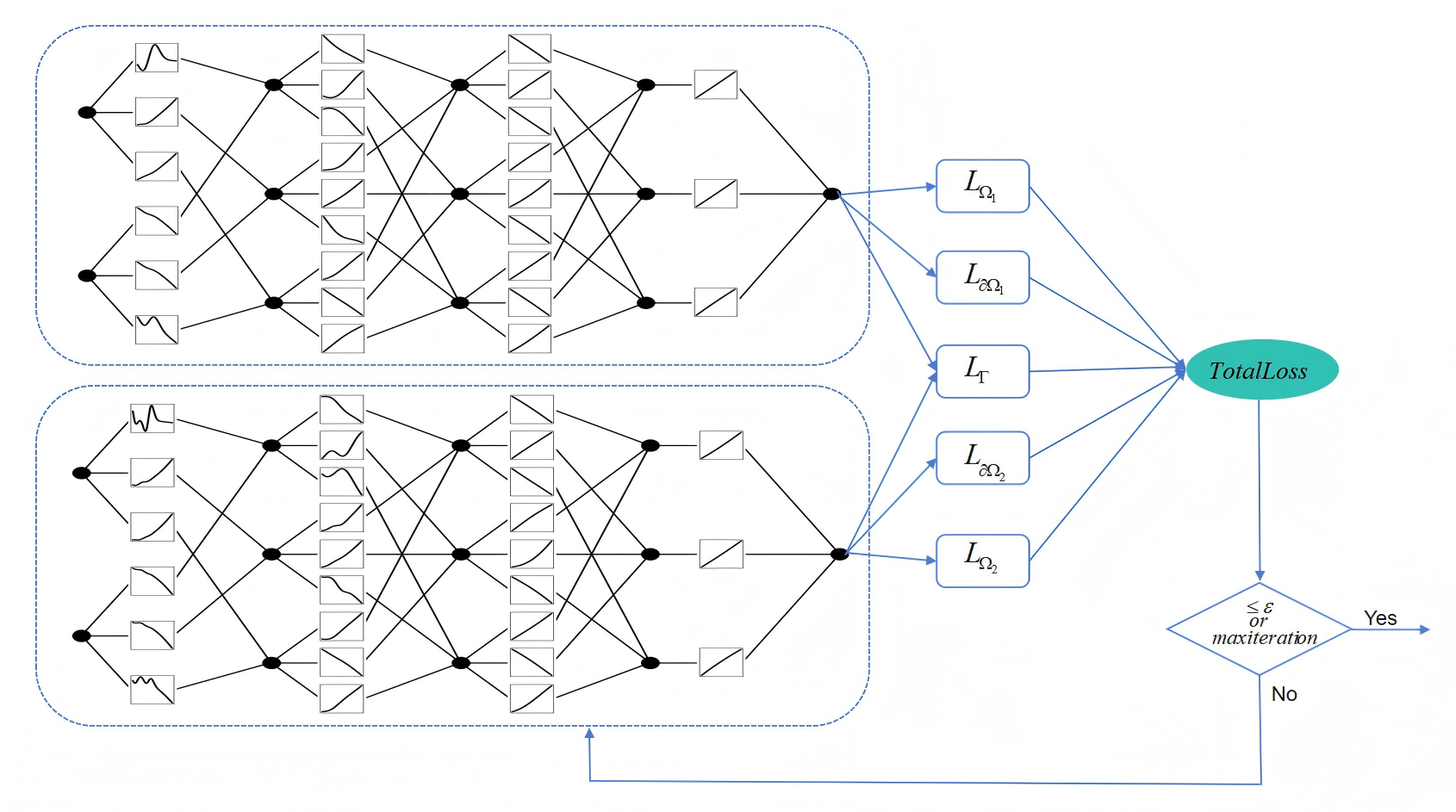} % 替换为图片的实际路径
    \caption{Schematic illustration of a dual KAN with structure [2,3,3,3,1].}
    \label{fig2:Schematic_KANs}
\end{figure}

In the original implementation of KANs, the univariate activation function is expressed as
\begin{align}\label{phi}
\phi(x)=c_{r} r(x)+c_{B} B(x),
\end{align}
here
\begin{align}
r(x)=\frac{x}{1+\exp (-x)}
\end{align}
is a residual-like activation function, and
\begin{align}
B\left(x\right)=\sum_{i=1}^{G+m}c_{i}B_{i}\left(x\right)
\end{align}
is an \( m \)-th order B-spline activation function defined on a grid with \( G \) intervals. For given \( G \) and \( m \), a set of spline basis functions \(\{B_i\}_{i=1}^{G+m}\) is uniquely defined, forming the activation function about Eq.\eqref{phi} and depending on the grid intervals \( G \). The parameters \( c_r, c_B \), and \(\sum_{i=1}^{G+m} c_i\) are trainable, thus the activation functions in KANs are not fixed, unlike in MLPs. This construction allows the activation functions to adapt to the specific features of the problem during training, providing a significant advantage over traditional fixed activations like ReLU or Tanh. A dual  KAN framework with structure [2,3,3,3,1] is schematically described in Figure \ref{fig2:Schematic_KANs}, for solving such interface problem.

Considering the discontinuity of the solution across the interface, we propose a dual KANs architecture to independently approximate the solution on both sides of the interface for solving elliptic interface problems. Specifically, we construct two separate KAN structures in $\Omega_{1}$ and $\Omega_{2}$ to approximate  $u|_{\Omega_1}$ and $u|_{\Omega_2}$ are constructed as $u_i(\mathbf{x}),\, \mathbf{x} \in \Omega_i, i=1, 2$, namely
\begin{align}
u_i(\mathbf{x})\approx\mathrm{KAN}_{i, \mathbf{\Phi}}(\mathbf{x})=(\mathbf{\Phi}_i^{L}\circ\mathbf{\Phi}_i^{L-1}\circ\cdots\circ\mathbf{\Phi}_i^{2}\circ\mathbf{\Phi}_i^{1})(\mathbf{x}), \quad i=1, 2.
\end{align}

In elliptic interface problems, the training residual is often highly non-uniform in space. This effect becomes pronounced in several common regimes, for example: (i) large contrast in coefficients $a_1,~a_2$, (ii) strong or spatially varying jump data $(\varphi,\psi)$, and (iii) curved interfaces where local solution features change rapidly along $\Gamma$. In such cases, uniform collocation may allocate too few points near the interface neighborhood where the residual is dominant, leading to slow convergence and interface-localized errors. This motivates the use of residual-driven adaptive sampling, and we adopt the RAR-D strategy described below.

\subsection{Residual-based adaptive refinement with diversity (RAR-D)}

The method presented in \cite{40} proposes RAR-D adaptive sampling  for improving the performance of Physics-Informed Neural Networks (PINNs). The method aim to enhance the accuracy of PINNs by dynamically adjusting the distribution of residual points based on the PDE residuals during training, which leads to more efficient training with fewer residual points. At the same time, this paper applies the RAR-D adaptive sampling method to the KAN network.

The concept of a density function $p(\mathbf{x})$  is designed to be proportional to the interface residual of the partial differential equation, and all residuals are resampled based on this density function. The newly generated residual points are thus distributed according to the density function. For any  point $\mathbf{x}$, we first compute the PDE residual $\varpi^{i}(\mathbf{x}) =  \left| -\nabla \cdot (a_{i} \nabla u_{i}(x^{j};\theta)) - f_{i} \right|$, and then compute a probability as
\begin{align}\label{midufun}
p(\mathbf{x}) \propto \frac{\varpi^{k}(\mathbf{x})}{\mathbb{E}\left[\varpi^{k}(\mathbf{x})\right]}+c. 
\end{align}
Here, \( k \geq 0 \) and \( c \geq 0 \) are two hyperparameters. $\mathbb{E}\left[\varpi^{k}(\mathbf{x})\right]$ can be approximated by a numerical integration such as Monte Carlo integration. The specific
adaptive sampling algorithm is demonstrated as Algorithm 1.
%
%\begin{algorithm}
%\renewcommand{\thealgorithm}{1} % 设置算法编号
%\caption{RAR-D\cite{40}}
%\begin{algorithmic}
%\State \textbf{Step1:} Sample the initial residual points \( \mathcal{T} \) using Latin hypercube sampling method;
%\State \textbf{Step2:} Train the network model for a certain number of iterations;
%\State \textbf{Step3:} Calculate the residuals of the initial interface residual points by automatically differentiating the partial differential equation;
%    \begin{enumerate}
%        \item Use the density function $p(\mathbf{x})$ from Eq.(16) to select new residual points, adding to the initial residual points.
%            \begin{enumerate}
%                \item Using the Latin hypercube sampling method, generate a series of point sets $\tau$;
%                \item Compute $p(\mathbf{x})$ for the points in $\tau$;
%                \item The normalized probability density function $p(\mathbf{x})$ defines a new probability density function $\tilde{p}(\mathbf{x})=\frac{p(\mathbf{x})}{A}$, with its normalization constant $A=\sum_{\mathbf{x} \in \tau} p(\mathbf{x})$;
%                \item According to $\tilde{p}(\mathbf{x})$, $m$ new residual points is selected from $\tau$.
%            \end{enumerate}
%    \end{enumerate}
%\State \textbf{Step5:} A certain number of iterations of network model training.
%\State \textbf{Step6:} the total number of iterations or the total number of residual points reaches the limit.
%\State \textbf{End for}
%\end{algorithmic}
%\end{algorithm}

\begin{algorithm}
\renewcommand{\thealgorithm}{1} % 设置算法编号
\caption{RAR-D \cite{40}}
\begin{algorithmic}
\State \textbf{Step 1:} Sample the initial sampling points \( \mathcal{T} \) using Latin hypercube sampling method;
\State \textbf{Step 2:} Train the network model for a certain number of iterations;
\State \textbf{Step 3:} Calculate the residuals of the initial sampling points by automatically differentiating the partial differential equation;
\State \textbf{Step 4:} Use the density function $p(\mathbf{x})$ from Eq. \eqref{midufun} to select new sampling points, replacing the initial residual points: 
            \begin{enumerate}
                \item Compute $p(\mathbf{x})$ for the points in \( \mathcal{T} \);
                \item The normalized probability density function $p(\mathbf{x})$ defines a new probability density function $\tilde{p}(\mathbf{x})=\frac{p(\mathbf{x})}{A}$, with its normalization constant $A=\sum_{\mathbf{x} \in \mathcal{T}} p(\mathbf{x})$;
                \item According to $\tilde{p}(\mathbf{x})$, the training set \( \mathcal{T} \) is updated.
            \end{enumerate}
\State \textbf{Step 5:} Training is continued until the total number of iterations reaches the limit or the residuals converge.
\State \textbf{End}
\end{algorithmic}
\end{algorithm}

\section{Numerical Experiments}

This section presents several numerical experiments for solving two-dimensional elliptic interface problems. First, numerical examples with analytical solutions were constructed to verify and compare the accuracy and effectiveness of PINNs and KANs. %in solving interface problems. 
By changing the number of neurons and \( G \) intervals of KANs, we investigated their influence on the solution of elliptic interface problems, thereby evaluating the impact of the KANs’ internal structure on the solution accuracy. Then, we constructed an analytical solution with discontinuities at the interface and compared the numerical results from two different network architectures. Then, an adaptive sampling algorithm was designed, and its integration with PINNs and KANs demonstrated the impact of adaptive sampling on solving elliptic interface problems.The results of these two models combined with the adaptive sampling algorithm are referred to as "PINNs-A" and "KANs-A", respectively.

To assess generalization, we evaluate the trained networks on an independent test set of $N_{test}$ in $\Omega$ randomly sampled points and report the maximum absolute error and the relative $\ell_2$-error, which are defined as
\begin{align}
|u-\hat{u}|_{\infty} = \max_{1 \leq j \leq N_{\text{test}}} |u(\mathbf{x}^j) - \hat{u}(\mathbf{x}^j)|, \quad
\|u-\hat{u}\|_{\ell_2} = \sqrt{\frac{\sum_{j=1}^{N_{\text{test}}} |u(\mathbf{x}^j) - \hat{u}(\mathbf{x}^j)|^2}{\sum_{j=1}^{N_{\text{test}}} |u(\mathbf{x}^j)|^2}}.
\end{align}
where $\hat{u}$ is the numerical solution function. In the numerical experiments hereafter, $e_{\Omega_{1}}$, $e_{\Omega_{2}}$, $e_{\Gamma}$, $e_{{\partial \Omega}_{1}}$ and  $e_{{\partial \Omega}_{2}}$ are the $\ell_2$-errors over \( \Omega_1 \), \( \Omega_2 \), interface \( \Gamma \) and boundary \( \partial \Omega_1 \), \( \partial \Omega_2 \), respectively.
\subsection{Continuous interface}
\begin{example}{\label{E1}}
%\textbf{Example 1.} 
As the first example, we consider the simple Poisson’s equation  in a square region $\Omega = [-1,1]\times[-1,1]$, with a circular interface $\Gamma=\{\mathbf{x}:x^{2}+y^{2}=0.25\}$, which cuts the domain into two sub-regions $\Omega_{1}=\{\mathbf{x}:\quad x^{2}+y^{2}\leq0.25\} $ and %$\Omega_2=\{\mathbf{x}:\quad x^2+y^2>0.25\}$
$\Omega_2=\Omega\setminus\Omega_1$. The coefficient is chosen to be $a_1= a_2=1$, and the analytical solution is selected as 
\begin{equation}
u(\mathbf{x}) = 
\begin{cases}
1 & \text{in } \Omega_1, \\
1 - \log_2 \sqrt{x^2 + y^2} & \text{in } \Omega_2.
\end{cases}
\end{equation}
The source terms $f_i$, boundary data $g_i$, and jump data $\psi,\varphi$ are obtained by substituting the above exact solution into \eqref{eq:model}.
\end{example}
Both PINNs and KANs are constructed to solve this problem. A dual PINNs with 3 layers and 20 neurons in each layer, and a dual KANs with 3 layers and 3 neurons in each layer are used. In KANs, we use $G=10$ to control the fineness of the spline function. Moreover, Latin hypercube sampling with $N_1=200$ in domain $ \Omega_{1}$, $N_2=500$ in domain $ \Omega_{2}$, $N_{ \Gamma}=300$ on the interface $\Gamma$, and $N_{\partial \Omega_{2}}=800$ on the boundary $\partial \Omega_{2}$, is applied. 

\begin{figure}[htbp]
    \centering
    
    % First row: single image with fixed height
    \begin{subfigure}[b]{0.4\textwidth}
        \centering
        \includegraphics[height=4cm]{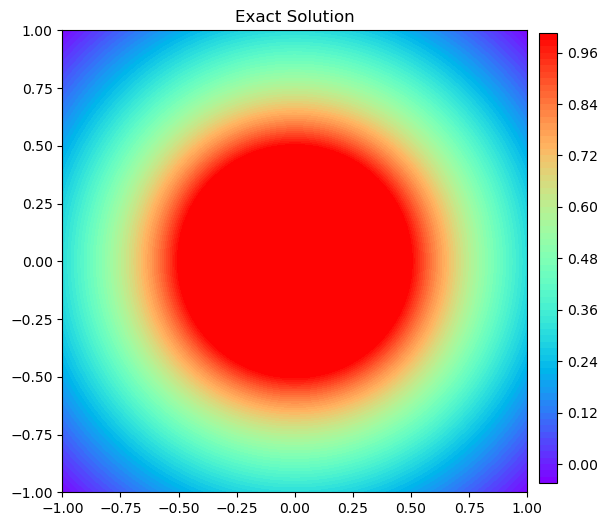} % Adjusted to use fixed height
        \caption{}
    \end{subfigure}

    % Adjust vertical spacing
    \vspace{0.2cm} 

    % Second row: two images with the same fixed height
    \begin{subfigure}[b]{0.4\textwidth}
        \centering
        \includegraphics[height=4cm]{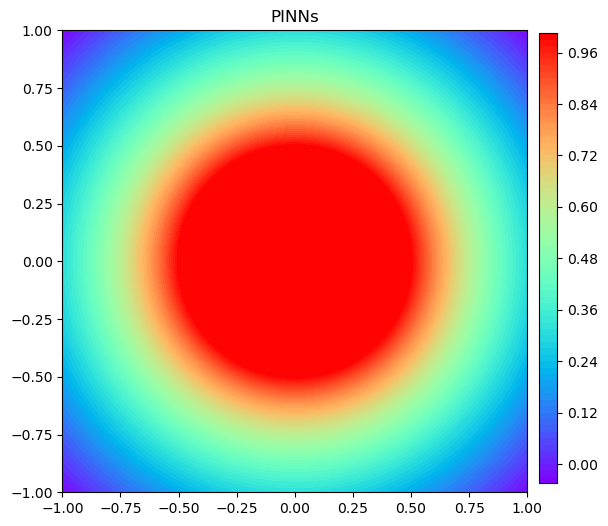} % Fixed height for consistency
        \caption{}
    \end{subfigure}
    \hspace{0.1cm} % Horizontal space between images
    \begin{subfigure}[b]{0.4\textwidth}
        \centering
        \includegraphics[height=4cm]{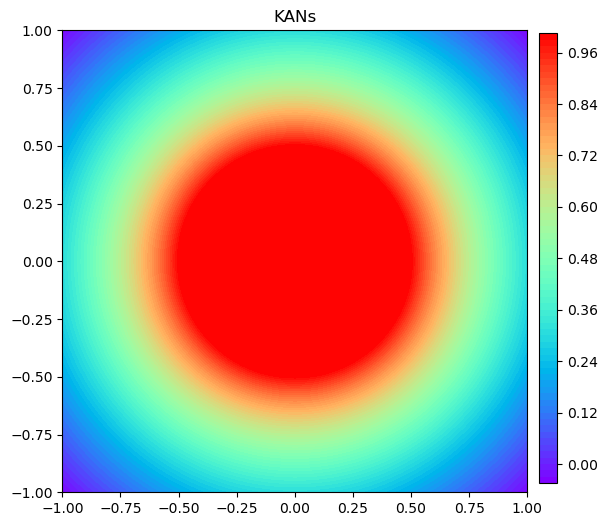} % Fixed height for consistency
        \caption{}
    \end{subfigure}

    % Adjust vertical spacing
    \vspace{0.2cm} 

    % Third row: three images with the same fixed height
    \begin{subfigure}[b]{0.4\textwidth}
        \centering
        \includegraphics[height=4cm]{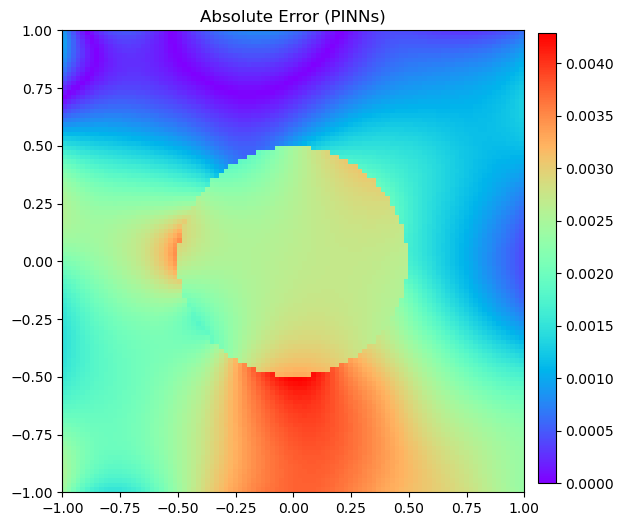} % Fixed height for consistency
        \caption{}
    \end{subfigure}
    \hspace{0.1cm} % Horizontal space between images
    \begin{subfigure}[b]{0.4\textwidth}
        \centering
        \includegraphics[height=4cm]{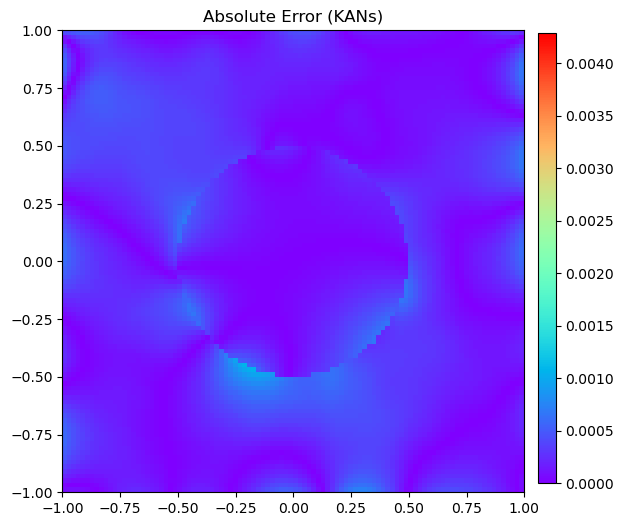} % Fixed height for consistency
        \caption{}
    \end{subfigure}

    \caption{Contour plots of the Example \ref{E1}: exact solution(a), approximation solutions by PINN(b) and KAN(c), and the absolute error by PINN(d) and KAN(e).}
    \label{fig:comparison3}
\end{figure}

\begin{figure}[htbp]
    \centering
    \includegraphics[width=0.45\textwidth]{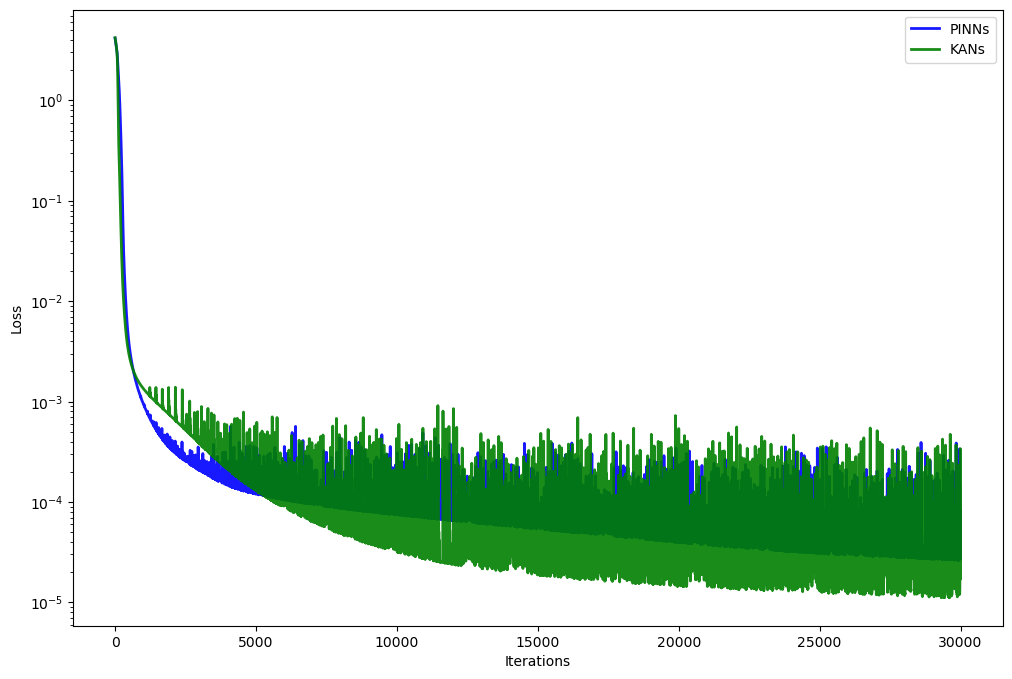} % 替换为图片的实际路径
    \caption{The evolution of loss error for Example \ref{E1}.}
    \label{fig:evolution4}
\end{figure}

Firstly, the contour plot of the analytical solution is shown in Fig.\ref{fig:comparison3}(a), and the approximate solutions obtained by PINN and KAN are given in Fig.\ref{fig:comparison3}(b)-(c), respectively. The corresponding absolute error solutions are plotted in Fig.\ref{fig:comparison3}(d)-(e), respectively. We can observe that both the proposed PINN and KAN can capture the interface by approximating the solution of the interface problem; meanwhile, KAN demonstrates significantly better approximation performance compared to PINN due to the error plots. We also plot the evolution of the loss function for both PINN and KAN solutions in Fig.\ref{fig:evolution4}, which indicates that KAN can accelerate the convergence of the loss function using a smaller network compared to PINN.

\begin{figure}[htbp]
    \centering
    \includegraphics[width=0.75\textwidth]{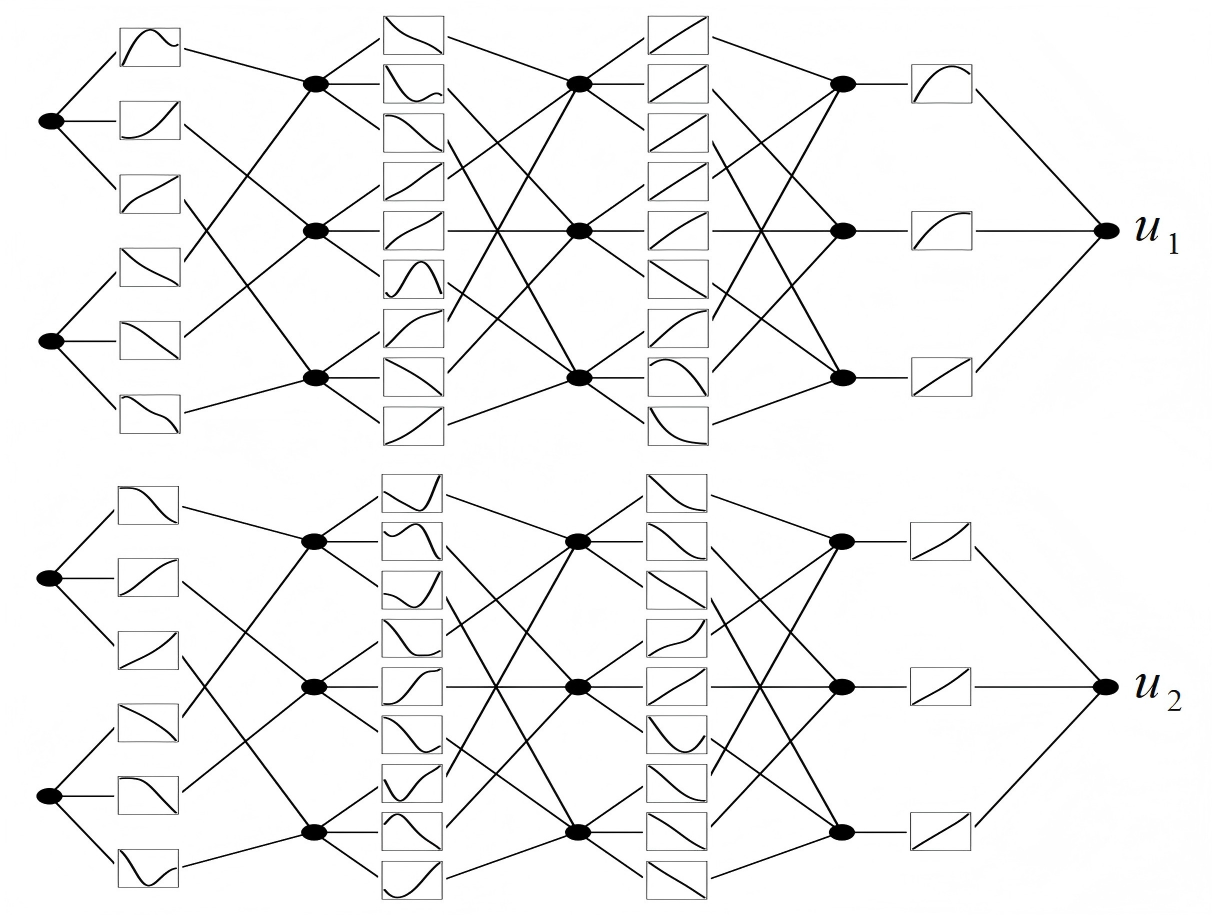} % 替换为图片的实际路径
    \caption{Visualization of the activation functions of a dual KAN after training for Example \ref{E1}.}
    \label{fig:KANexample1}
\end{figure}

Fig.\ref{fig:KANexample1}  illustrates the activation functions of the coupled KAN after training, while the initial activation functions used here can check Fig.\ref{fig2:Schematic_KANs} by comparison. It can be observed that the basis functions in each layer undergo adaptive adjustments during the training process, thereby forming diverse nonlinear representation capabilities. This mechanism ensures that the network maintains strong expressive power and convergence efficiency while keeping a relatively small scale.

The performance of KAN depends on the network architecture and the grid intervals \( G \). To this end, Table \ref{table:t1} shows the relative $\ell_2$ error with \( G=10 \) for different numbers of neurons, and
Table \ref{table:t2} lists the result by fixing neuron=5 and different numbers of \( G \) intervals in KAN. The results demonstrate that neither the number of neurons nor \( G \) intervals significantly impacts computational error when one parameter is fixed while varying the other.
 %Table \ref{table:t2} shows the relative $L_2$ error with fixed neuron=5 for different numbers of grids in KAN. From the table, it can be seen that increasing the number of neurons and grid size does not necessarily lead to better computational accuracy. 
 Therefore, an appropriate network depth and width can be selected to achieve good accuracy and efficiency. 
 
During the adaptive sampling step, both networks of PINN and KAN are trained using the Adam for 40000 steps.  We generate 5000 initial sampling points for training 20000 steps, subsequently adaptive resampling (RAR-D) is employed after every 2000 steps. Here, $k$ = 2 and $c$= 0.  Under the same initial sampling, we apply the PINN and KAN trainings, and the RAR-D adaptive algorithm, then we display the computational results in
Table  \ref{table:t3}. From this table, we can conclude that the adaptive sampling method of both networks effectively improves the computational accuracy in the $\ell_2$-norm. More importantly, we can see that each error in the second row for KAN is much smaller than corresponding error in the first row for PINN. Also in Fig. \ref{fig:comparison3}(d)-(e), the computational error distribution by KAN is much more uniform than that of PINN.
These observations clearly demonstrate that the performance of KAN in solving the elliptic interface problem is superior to that of PINN.

%Table  \ref{table:t3} shows the relative $L_2$ error of the same initial sampling training and adaptive sampling training. It is shown that the adaptive sampling method effectively reduces the $L_2$ error. The results from this table clearly demonstrate that the performance of KANs in solving the problem is superior to that of PINNs.

\begin{table}[H]
\centering
\caption{Comparison of errors in different numbers of neurons in KANs for Example \ref{E1}.}\label{table:t1}
\begin{tabular}{c c c c c}
\toprule
\textbf{Number of Neurons} & $\boldsymbol{e_{\Omega_{1}}}$ & $\boldsymbol{e_{\Omega_{2}}}$ & $\boldsymbol{e_{\Gamma}}$ & $\boldsymbol{e_{\partial \Omega_{2}}}$  \\
\midrule
3 &  $1.456 \times 10^{-4}$ & $5.194 \times 10^{-4}$ & $3.402 \times 10^{-4}$ & $3.022 \times 10^{-3}$ \\
5 & $2.565 \times 10^{-4}$ & $3.982 \times 10^{-4}$ & $3.064 \times 10^{-4}$ & $1.124 \times 10^{-3}$  \\
7 & $2.074 \times 10^{-4}$ & $7.310 \times 10^{-4}$ & $2.587 \times 10^{-4}$ & $1.485 \times 10^{-3}$  \\
\bottomrule
\end{tabular}
\end{table}

\begin{table}[H]
\centering
\caption{Comparison of errors in different \( G \) intervals in KANs for Example \ref{E1}.}\label{table:t2}
\begin{tabular}{c c c c c}
\toprule
\textbf{\( G \) intervals} & $\boldsymbol{e_{\Omega_{1}}}$ & $\boldsymbol{e_{\Omega_{2}}}$ & $\boldsymbol{e_{\Gamma}}$ & $\boldsymbol{e_{\partial \Omega_{2}}}$  \\
\midrule
5 &  $5.232 \times 10^{-4}$ & $9.185 \times 10^{-4}$ & $1.002 \times 10^{-3}$ & $2.421 \times 10^{-3}$ \\
10 & $2.565 \times 10^{-4}$ & $3.982 \times 10^{-4}$ & $3.064 \times 10^{-4}$ & $1.124 \times 10^{-3}$  \\
15 & $4.298 \times 10^{-5}$ & $3.281 \times 10^{-4}$ & $2.343 \times 10^{-4}$ & $1.677 \times 10^{-3}$ \\
\bottomrule
\end{tabular}
\end{table}

\begin{table}[H]
\centering
\caption{Comparison of errors between training with the same initial sampling and RAR-D adaptive sampling for Example \ref{E1}.}\label{table:t3}
\begin{tabular}{c c c c c c}
\toprule
\textbf{Networks} & $\boldsymbol{e_{\Omega_{1}}}$ & $\boldsymbol{e_{\Omega_{2}}}$ & $\boldsymbol{e_{\Gamma}}$ & $\boldsymbol{e_{\partial \Omega_{2}}}$ & $|u-\hat{u}|_{\infty}$ \\
\midrule
PINNs  & $9.790 \times 10^{-4}$ & $3.280 \times 10^{-3}$ & $1.046 \times 10^{-3}$ &  $8.297 \times 10^{-3}$ &  $3.693 \times 10^{-3}$ \\
KANs & $2.040 \times 10^{-4}$ & $5.024 \times 10^{-4}$ & $2.889 \times 10^{-4}$ & $1.156 \times 10^{-3}$ &  $8.079 \times 10^{-4}$  \\
PINNs-A & $5.634 \times 10^{-4}$ & $8.883 \times 10^{-4}$ & $6.641 \times 10^{-4}$ & $2.842 \times 10^{-3}$ &  $1.075 \times 10^{-3}$  \\
KANs-A & $1.035 \times 10^{-4}$ & $1.927 \times 10^{-4}$ & $1.680 \times 10^{-4}$ & $6.512 \times 10^{-4}$ &   $3.437 \times 10^{-4}$ \\

\bottomrule
\end{tabular}
\end{table}

\subsection{Discontinuous interface}

\begin{figure}[htbp]
    \centering
    
    % First row: single image with fixed height
    \begin{subfigure}[b]{0.4\textwidth}
        \centering
        \includegraphics[height=4cm]{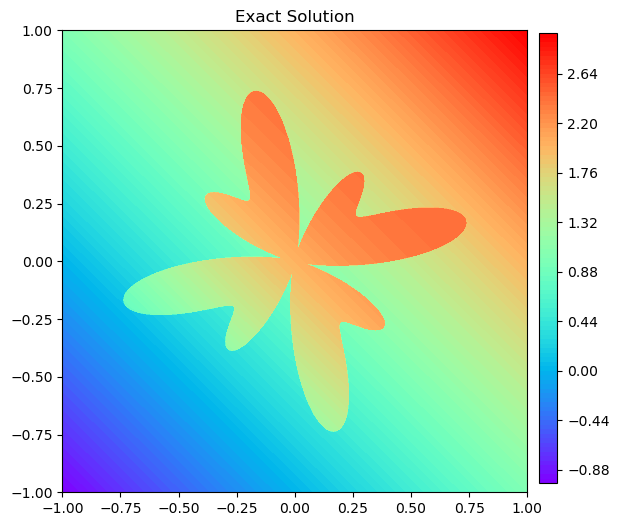} % Adjusted to use fixed height
        \caption{}
    \end{subfigure}

    % Adjust vertical spacing
    \vspace{0.2cm} 

    % Second row: two images with the same fixed height
    \begin{subfigure}[b]{0.4\textwidth}
        \centering
        \includegraphics[height=4cm]{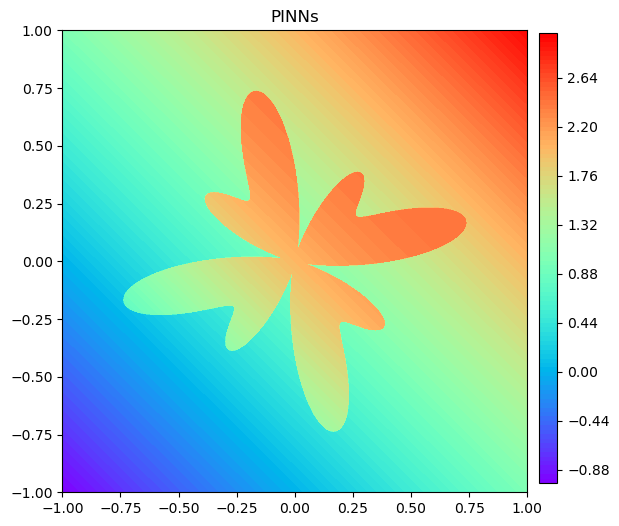} % Fixed height for consistency
        \caption{}
    \end{subfigure}
    \hspace{0.1cm} % Horizontal space between images
    \begin{subfigure}[b]{0.4\textwidth}
        \centering
        \includegraphics[height=4cm]{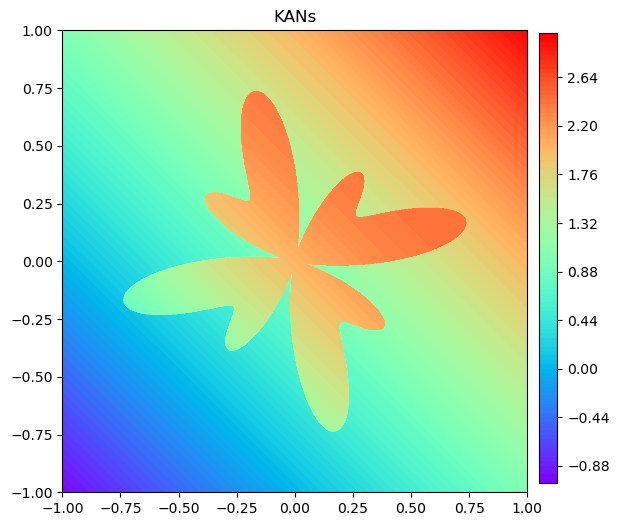} % Fixed height for consistency
        \caption{}
    \end{subfigure}

    % Adjust vertical spacing
    \vspace{0.2cm} 

    % Third row: three images with the same fixed height
    \begin{subfigure}[b]{0.4\textwidth}
        \centering
        \includegraphics[height=4cm]{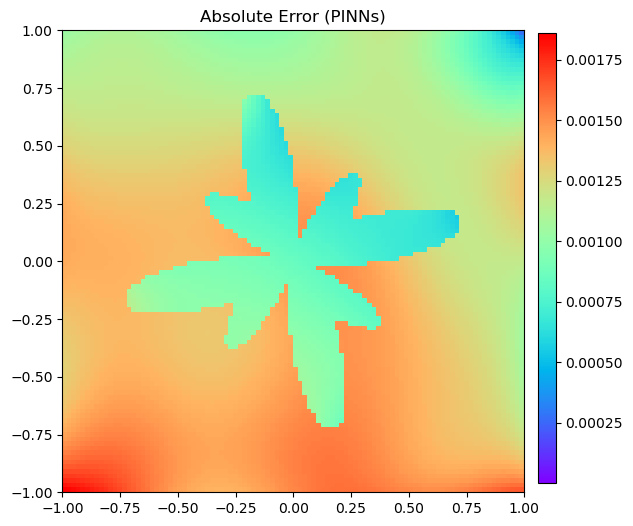} % Fixed height for consistency
        \caption{}
    \end{subfigure}
    \hspace{0.1cm} % Horizontal space between images
    \begin{subfigure}[b]{0.4\textwidth}
        \centering
        \includegraphics[height=4cm]{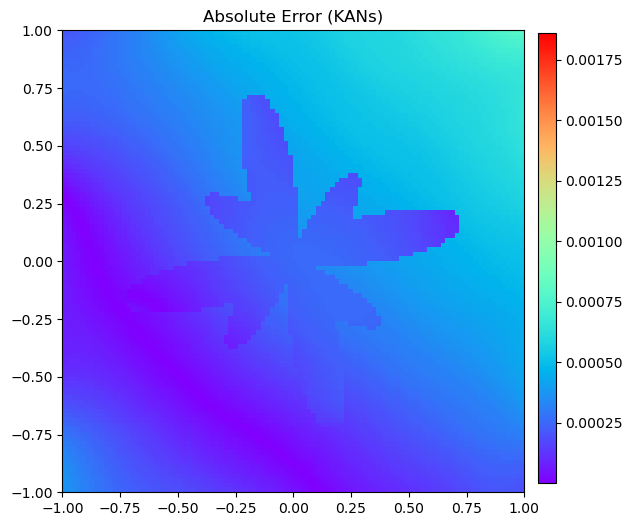} % Fixed height for consistency
        \caption{}
    \end{subfigure}

    \caption{Contour plots of the Example \ref{E4}, exact solution(a), approximation solutions by PINNs(b) and KANs(c), and the absolute error by PINNs(d) and KANs(e).}
    \label{fig:comparison9}
\end{figure}

\begin{figure}[htbp]
    \centering
    \includegraphics[width=0.45\textwidth]{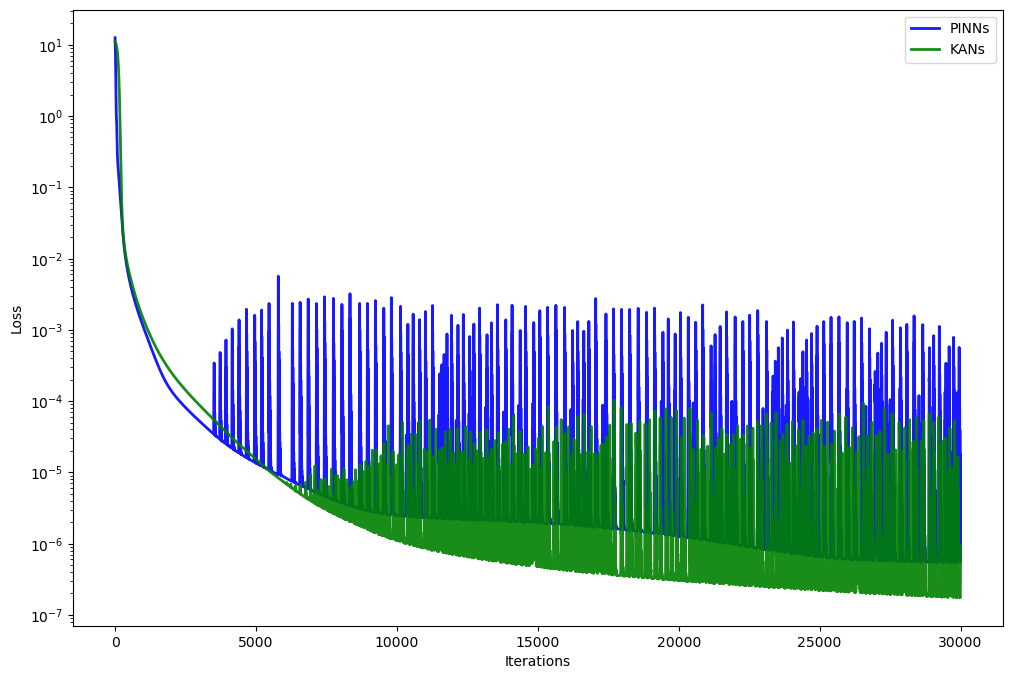} % 替换为图片的实际路径
    \caption{The evolution of loss error for Example \ref{E4}.}
    \label{fig:evolution10}
\end{figure}

\begin{example}{\label{E4}}
%\textbf{Example 4.} 
This example employs a more sophisticated interface, parameterized by
\begin{equation}
\begin{cases}
x(\theta) = (a + b \cos(m \theta)) \sin(n \theta) \cos(\theta), \\
y(\theta) = (a + b \cos(m \theta)) \sin(n \theta) \sin(\theta),
\end{cases}
\end{equation}
where $\theta\in \left [ 0,2\pi  \right ] $, $a=b=0.40178$, $m=2$ and $n=6$. In the domain $\Omega = [-1,1]\times[-1,1]$, $\Omega_{1}$ and $\Omega_{2}$ are defined,  as the region inside and outside the interface $\Gamma$, respectively. The coefficient $a_i$ and the solution are given as follows
\begin{equation*}
a = \left\{
\begin{array}{ll}
\frac{x^{2} - y^{2} + 3}{7},  & \text { in } \Omega_{1}, \\
\frac{2 + xy}{5},  & \text { in } \Omega_{2}, \\  
\end{array}
\right.
\quad
u(\mathbf{x})= \left\{
\begin{array}{ll}
\sin(x + y) + \cos(x + y) + 1,  & \text { in } \Omega_{1}, \\
x + y + 1,  & \text { in } \Omega_{2}.
\end{array}
\right.
\end{equation*}
\end{example}
A contour plot of the analytical solution of this problem is plotted in Fig.\ref{fig:comparison9}(a). To solve this Example \ref{E4}, a dual PINN with 3 layers and 20 neurons in each layer and a dual KAN with 3 layers and 3 neurons in each layer are constructed. 
Select $N_1=300,N_2=500,N_{ \Gamma}=300,N_{\partial \Omega_{2}}=800$ for both networks, and choose $G=5$ in KAN.
Similar with Example \ref{E1}, we show the contour plots of the PINN and KAN solutions in Fig.\ref{fig:comparison9}(b)-(c), and the corresponding absolute error in  \ref{fig:comparison9}(d)-(e). Also, Fig.\ref{fig:evolution10} presents the evolution of the loss function of PINN and KAN, which indicates faster convergence of KAN compared with PINN. 

Tables \ref{table:t10} and \ref{table:t11} list  the performance of KAN along with the change of the numbers of neurons or $G$ intervals. From these two tables, in the sense of $L^2$-norm error, 3 neurons are suitable for KAN with $G=5$, and such KAN is much less sensitive with the scale of  $G$ intervals, so $G=5$ is enough under present case.

As for the adaptive sampling, both networks of PINN and KAN are trained using the Adam for 30000 steps. During the Adam training epoch, we first generates 5000 residual points for training 20000 steps, subsequently applying adaptive resampling (RAR-D) after every 1000 iterations. Here, $k$ = 2 and $c$= 1. Table \ref{table:t12} shows the comparisons of PINN and KAN, and their adaptive performances, in the sense of $L^2$-error. 

Evidently, the adaptive algorithm demonstrates improved performance for both networks, with KAN outperforming PINN under the same initial sampling conditions.
%Slightly different with the previous examples, a dual PINN with 3 layers and 20 neurons (with more numbers), and a dual KAN with 3 layers and 5 neurons in each layer are constructed. Actually, as indicated in both Examples \ref{E1} and \ref{E2}, combination of error results in Table \ref{table:t10} (with $G=10$) and \ref{table:t11} (under different number of $G$) recommends the selection of 3-5 neurons for KAN, and $G=5$. The contour plots of the solutions and evolution of loss error for both networks are presented in Fig. \ref{fig:comparison9} and \ref{fig:evolution10}. {Exactly same with previous Example?\color{red} In the adaptive sampling setup, We trained the neural network using the Adam  for 30,000 steps. We generate 5000 residual points and resampled every 1000 iterations for adaptive sampling after training 20000 steps.} 
%Here, $k$ = 2 and $c$= 1. Table \ref{table:t12} shows the relative $L_2$ error of the same initial sampling training and adaptive sampling training. It is shown that the adaptive sampling method effectively reduces the $L_2$ error. The results from this table clearly demonstrate that the performance of KANs in solving the problem is superior to that of PINNs.

\begin{table}[H]
\centering
\caption{Comparison of errors in different numbers of neurons in KANs for Example \ref{E4}}\label{table:t10}
\begin{tabular}{c c c c c}
\toprule
\textbf{Number of Neurons} & $\boldsymbol{e_{\Omega_{1}}}$ & $\boldsymbol{e_{\Omega_{2}}}$ & $\boldsymbol{e_{\Gamma}}$ & $\boldsymbol{e_{\partial \Omega_{2}}}$  \\
\midrule
3 &  $6.865 \times 10^{-5}$ & $5.646 \times 10^{-5}$ & $7.124 \times 10^{-5}$ & $8.841 \times 10^{-5}$ \\
5 & $2.779 \times 10^{-4}$ & $2.832 \times 10^{-4}$ & $2.714 \times 10^{-4}$ & $2.714 \times 10^{-4}$  \\
7 & $1.051 \times 10^{-4}$ & $8.843 \times 10^{-5}$ & $1.022 \times 10^{-4}$ & $7.650 \times 10^{-5}$  \\
\bottomrule
\end{tabular}
\end{table}

\begin{table}[H]
\centering
\caption{Comparison of errors in different \( G \) intervals in KANs for Example \ref{E4}.}\label{table:t11}
\begin{tabular}{c c c c c}
\toprule
\textbf{\( G \) intervals} & $\boldsymbol{e_{\Omega_{1}}}$ & $\boldsymbol{e_{\Omega_{2}}}$ & $\boldsymbol{e_{\Gamma}}$ & $\boldsymbol{e_{\partial \Omega_{2}}}$  \\
\midrule
5 &  $3.459 \times 10^{-5}$ & $2.892 \times 10^{-5}$ & $3.028 \times 10^{-5}$ & $3.453 \times 10^{-5}$ \\
10 & $2.779 \times 10^{-4}$ & $2.832 \times 10^{-5}$ & $2.714 \times 10^{-4}$ & $2.714 \times 10^{-5}$  \\
15 & $9.490 \times 10^{-5}$ & $1.031 \times 10^{-4}$ & $8.375 \times 10^{-5}$ & $1.236 \times 10^{-4}$  \\
\bottomrule
\end{tabular}
\end{table}

\begin{table}[H]
\centering
\caption{Comparison of errors between training with the same initial sampling and RAR-D adaptive sampling for Example \ref{E4}.}\label{table:t12}
\begin{tabular}{c c c c c c}
\toprule
\textbf{Networks} & $\boldsymbol{e_{\Omega_{1}}}$ & $\boldsymbol{e_{\Omega_{2}}}$ & $\boldsymbol{e_{\Gamma}}$ & $\boldsymbol{e_{\partial \Omega_{2}}}$ & $|u-\hat{u}|_{\infty}$ \\
\midrule
PINNs  & $2.935 \times 10^{-4}$ & $4.671 \times 10^{-4}$ & $2.767 \times 10^{-4}$ &  $4.188 \times 10^{-4}$ & $1.077 \times 10^{-3}$ \\

KANs & $5.357 \times 10^{-5}$ & $8.319 \times 10^{-5}$ & $5.255 \times 10^{-5}$ & $7.743 \times 10^{-5}$ & $2.251 \times 10^{-4}$ \\
PINNs-A & $1.012 \times 10^{-4}$ & $1.022 \times 10^{-4}$ & $1.031 \times 10^{-4}$ & $9.879 \times 10^{-5}$ & $4.225 \times 10^{-4}$\\
KANs-A & $3.648 \times 10^{-5}$ & $2.177 \times 10^{-5}$ & $3.679 \times 10^{-5}$ & $3.146 \times 10^{-5}$ & $1.512 \times 10^{-4}$  \\
\bottomrule
\end{tabular}
\end{table}

In the following, we present two more examples of elliptic interface problems with complex interface/boundary geometries.
\begin{example}{\label{E5}}
The computational domain consists of two highly irregular, non-convex subdomains, with the interface given by
\begin{align}
r=0.6+0.216\sin(3\theta)+0.096\cos(2\theta)+0.24\cos(5\theta),
\end{align}
The boundary is given in polar coordinates
\begin{align}
r=0.14\sin(4\theta)+0.12\cos(6\theta)+0.09\cos(5\theta),
\end{align}
where $\theta\in \left [ 0,2\pi  \right ] $. Over the domain $\Omega = [-2,2]\times[-2,2]$, $\Omega_{1}$ and $\Omega_{2}$ are defined, respectively, to be the region inside and outside $\Gamma$. The coefficients $a_i$ and the solutions are given as follows
\begin{equation}
a = \left\{
\begin{array}{ll}
x^{2} + y^{2},  & \text { in } \Omega_{1}, \\
 xy,  & \text { in } \Omega_{2}, \\  
\end{array}
\right.
\quad
u(\mathbf{x})= \left\{
\begin{array}{ll}
\cos(x + y) ,  & \text { in } \Omega_{1}, \\
\sin(x + y) ,  & \text { in } \Omega_{2}.
\end{array}
\right.
\end{equation}
\end{example}

\begin{example}{\label{E6}}
    Consider an annular region with inner and outer radii $r_{in}=0.151$ and $r_{out}=0.911$, and the immersed star-shaped interface $\gamma$ is described by the level-set function:
\begin{align}
\phi(x,y)=\sqrt{x^2+y^2}-r_0\left(1+\sum_{k=1}^3\beta_k\cos\left(n_k\left(\arctan\left(\frac{y}{x}\right)-\theta_k\right)\right)\right),
\end{align}
and
\begin{align*}
r_0=0.483,\quad\begin{pmatrix}n_1\\\beta_1\\\theta_1\end{pmatrix}=\begin{pmatrix}3\\0.1\\0.5\end{pmatrix},\quad\begin{pmatrix}n_2\\\beta_2\\\theta_2\end{pmatrix}=\begin{pmatrix}4\\-0.1\\1.8\end{pmatrix}\quad\text{and}\quad\begin{pmatrix}n_3\\\beta_3\\\theta_3\end{pmatrix}=\begin{pmatrix}7\\0.15\\0\end{pmatrix}.
\end{align*}
where $\theta\in \left [ 0,2\pi  \right ] $. The domain $\Omega = [-1,1]\times[-1,1]$ is decomposed by this interface into $\Omega_{1}$ and $\Omega_{2}$, where the coefficients $a_i$ and the solution are defined by
\begin{align}
a &= \left\{
\begin{array}{ll}
10\left(1+\frac{1}{5}\cos(2\pi(x+y))\sin(2\pi(x-y))\right),  & \text { in } \Omega_{1}, \\
1,  & \text { in } \Omega_{2},
\end{array}
\right. \\
u(\mathbf{x}) &= \left\{
\begin{array}{ll}
\sin(2x)\cos(2y),  & \text { in } \Omega_{1}, \\
\left(16\left(\frac{y-x}{3}\right)^5-20\left(\frac{y-x}{3}\right)^3+5\left(\frac{y-x}{3}\right)\right)\log\left(x+y+3\right),  & \text { in } \Omega_{2}.
\end{array}
\right.
\end{align}
\end{example}

\begin{figure}[htbp]
    \centering
    
    % First row: single image with fixed height
    \begin{subfigure}[b]{0.4\textwidth}
        \centering
        \includegraphics[height=4cm]{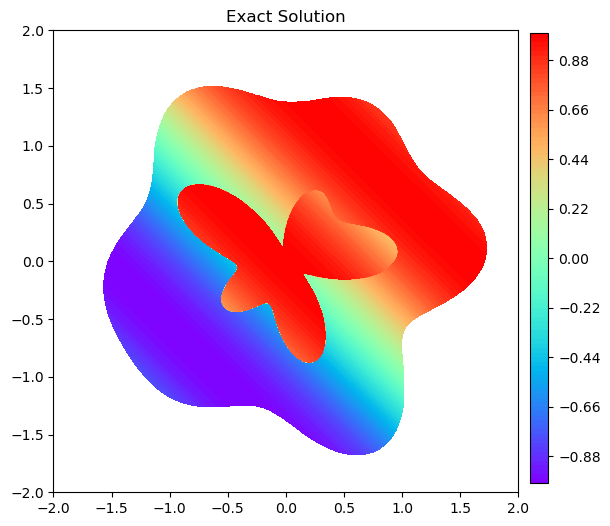} % Adjusted to use fixed height
        \caption{}
    \end{subfigure}

    % Adjust vertical spacing
    \vspace{0.2cm} 

    % Second row: two images with the same fixed height
    \begin{subfigure}[b]{0.4\textwidth}
        \centering
        \includegraphics[height=4cm]{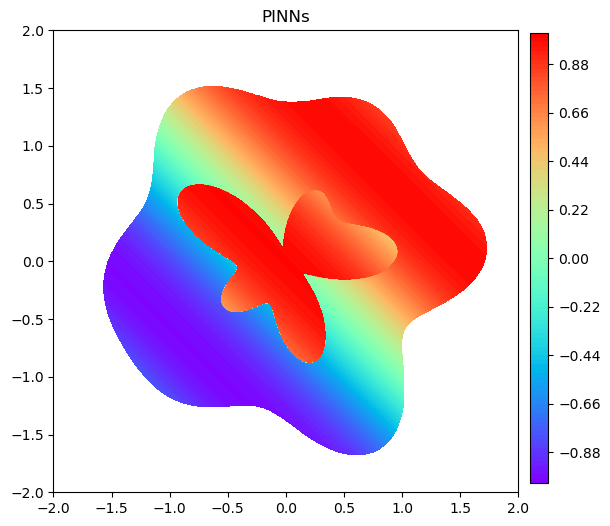} % Fixed height for consistency
        \caption{}
    \end{subfigure}
    \hspace{0.1cm} % Horizontal space between images
    \begin{subfigure}[b]{0.4\textwidth}
        \centering
        \includegraphics[height=4cm]{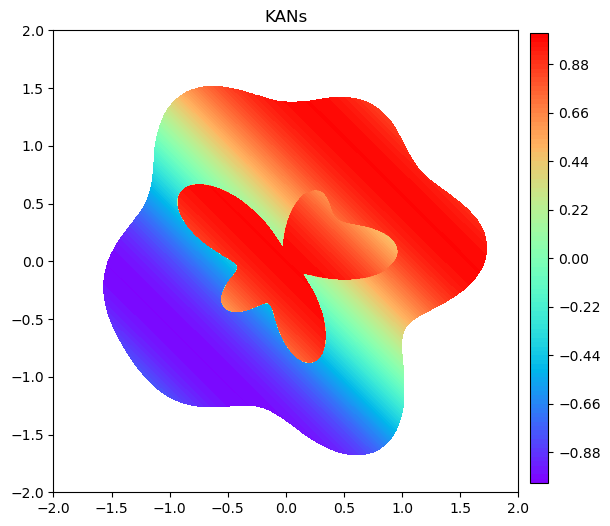} % Fixed height for consistency
        \caption{}
    \end{subfigure}

    % Adjust vertical spacing
    \vspace{0.2cm} 

    % Third row: three images with the same fixed height
    \begin{subfigure}[b]{0.4\textwidth}
        \centering
        \includegraphics[height=4cm]{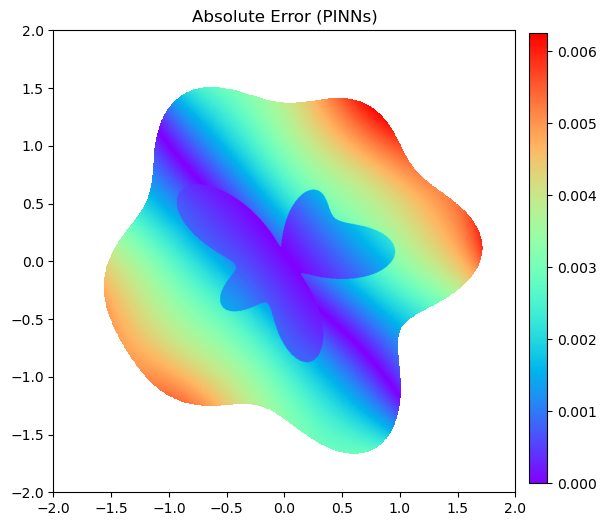} % Fixed height for consistency
        \caption{}
    \end{subfigure}
    \hspace{0.1cm} % Horizontal space between images
    \begin{subfigure}[b]{0.4\textwidth}
        \centering
        \includegraphics[height=4cm]{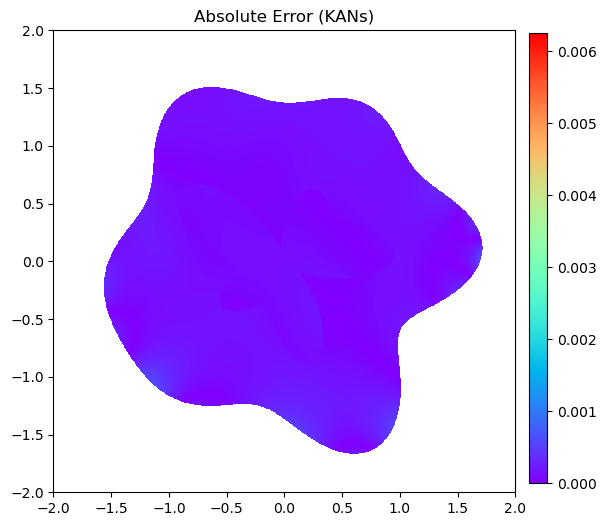} % Fixed height for consistency
        \caption{}
    \end{subfigure}

    \caption{Contour plots of the Example \ref{E5}, exact solution(a), approximation solutions by PINNs(b) and KANs(c), and the absolute error by PINNs(d) and KANs(e).}
    \label{fig:comparison11}
\end{figure}

\begin{figure}[htbp]
    \centering
    \includegraphics[width=0.45\textwidth]{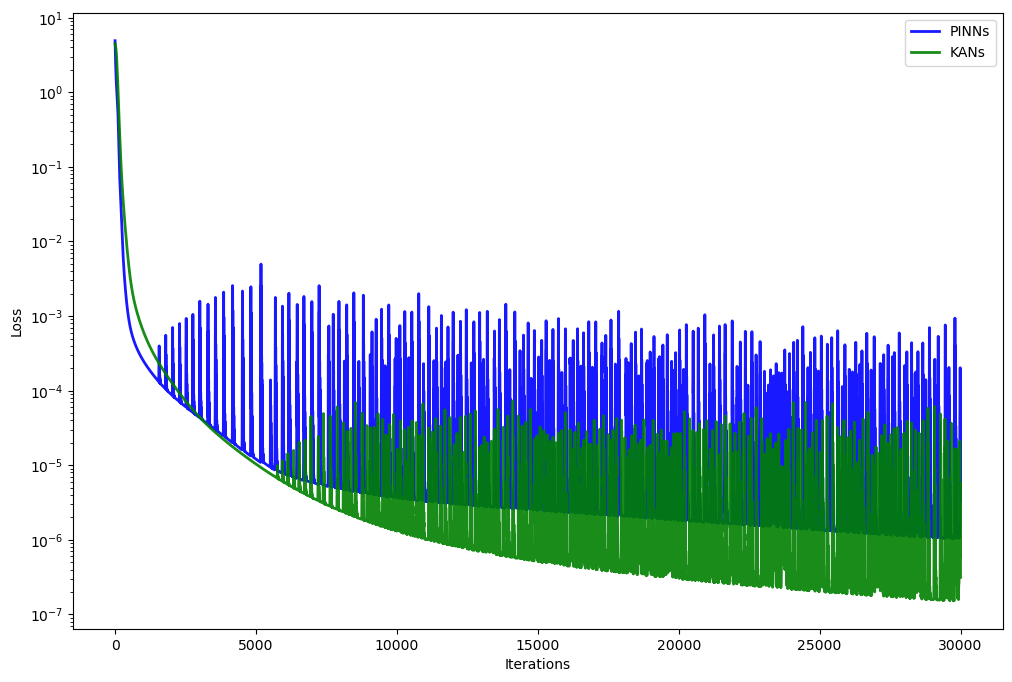} % 替换为图片的实际路径
    \caption{The evolution of loss error for Example \ref{E5}.}
    \label{fig:evolution12}
\end{figure}

\begin{figure}[htbp]
    \centering
    
    % First row: single image with fixed height
    \begin{subfigure}[b]{0.4\textwidth}
        \centering
        \includegraphics[height=4cm]{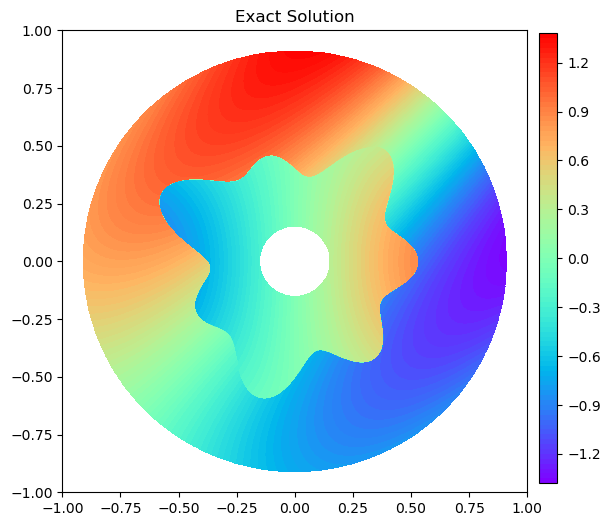} % Adjusted to use fixed height
        \caption{}
    \end{subfigure}

    % Adjust vertical spacing
    \vspace{0.2cm} 

    % Second row: two images with the same fixed height
    \begin{subfigure}[b]{0.4\textwidth}
        \centering
        \includegraphics[height=4cm]{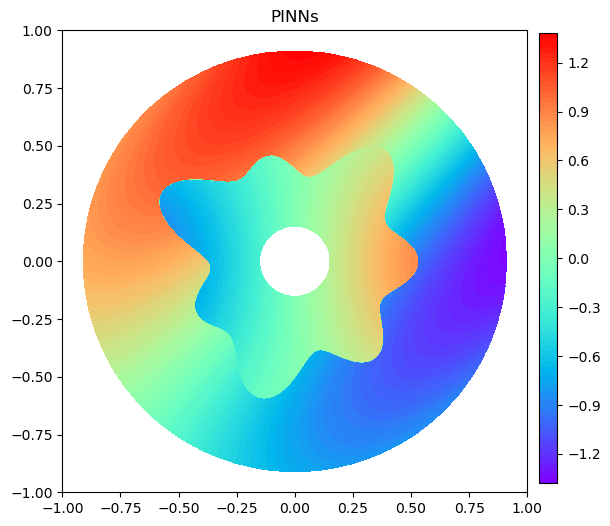} % Fixed height for consistency
        \caption{}
    \end{subfigure}
    \hspace{0.1cm} % Horizontal space between images
    \begin{subfigure}[b]{0.4\textwidth}
        \centering
        \includegraphics[height=4cm]{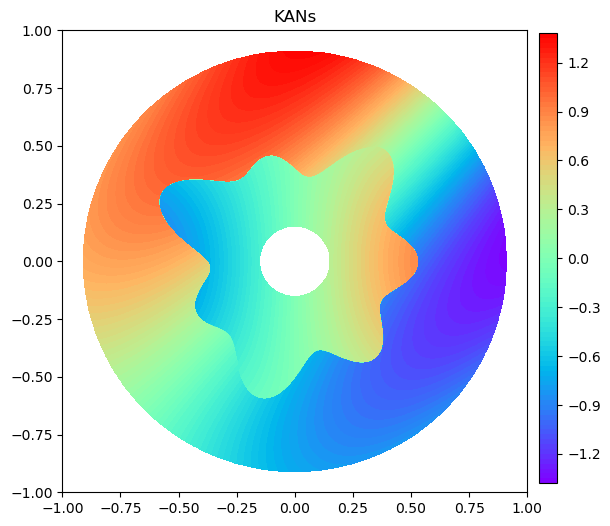} % Fixed height for consistency
        \caption{}
    \end{subfigure}

    % Adjust vertical spacing
    \vspace{0.2cm} 

    % Third row: three images with the same fixed height
    \begin{subfigure}[b]{0.4\textwidth}
        \centering
        \includegraphics[height=4cm]{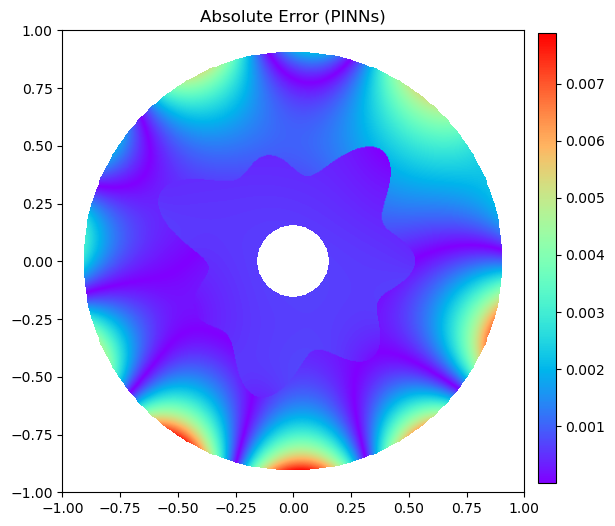} % Fixed height for consistency
        \caption{}
    \end{subfigure}
    \hspace{0.1cm} % Horizontal space between images
    \begin{subfigure}[b]{0.4\textwidth}
        \centering
        \includegraphics[height=4cm]{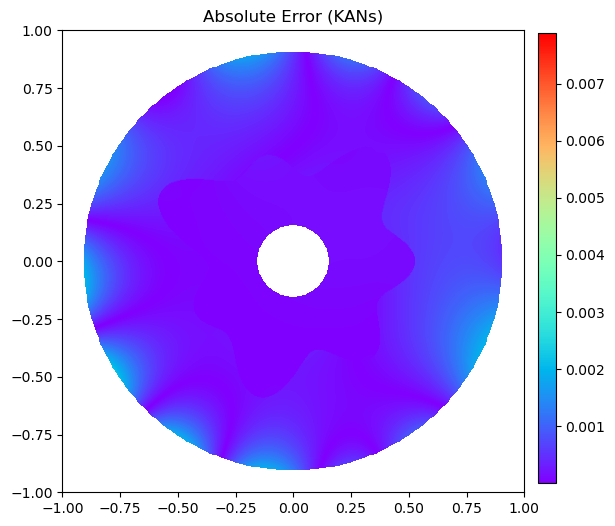} % Fixed height for consistency
        \caption{}
    \end{subfigure}

    \caption{Contour plots of the Example \ref{E6}, exact solution(a), approximation solutions by PINNs(b) and KANs(c), and the absolute error by PINNs(d) and KANs(e).}
    \label{fig:comparison13}
\end{figure}

\begin{figure}[htbp]
    \centering
    \includegraphics[width=0.45\textwidth]{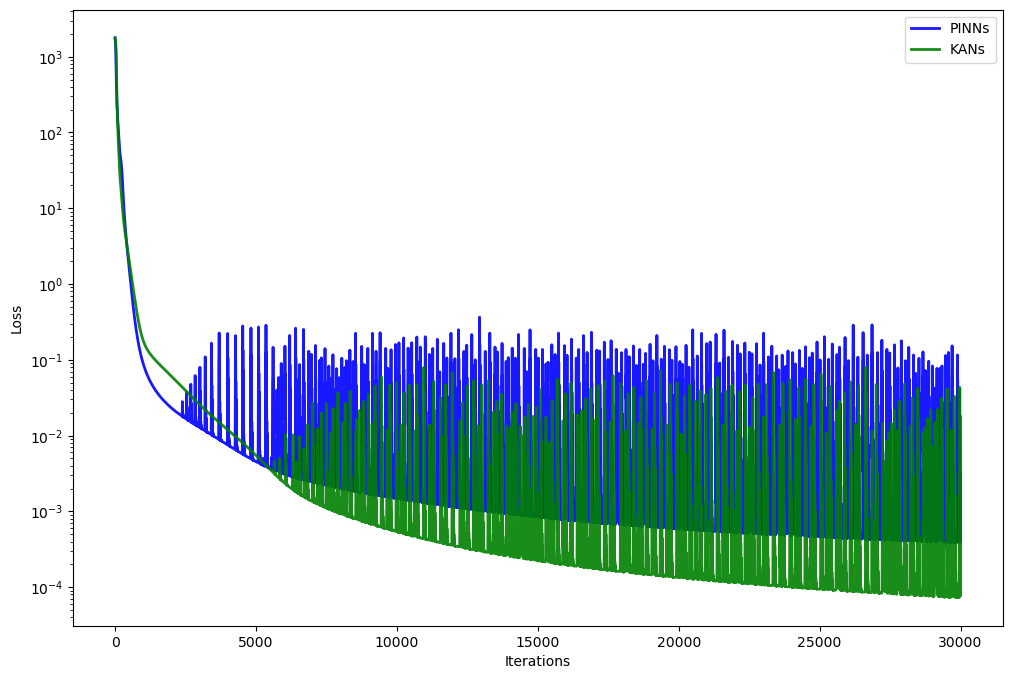} % 替换为图片的实际路径
    \caption{The evolution of loss error for Example \ref{E6}.}
    \label{fig:evolution14}
\end{figure}

For both examples, we use the same network settings, namely a dual PINN with 3 layers and 20 neurons in each layer and a dual KAN with 3 layers and 5 neurons in each layer. Contour plots of the analytical, numerical/error solutions and evolutions of the loss error are shown in Fig.\ref{fig:comparison11}- \ref{fig:evolution14}, respectively for Examples \ref{E5} and \ref{E6}.
Select $N_1=300,N_2=500,N_{ \Gamma}=300$, and $N_{\partial \Omega_{2}}=800$(for Example \ref{E5}), 300 (for Example \ref{E6}) for both networks, and choose $G=5$ in KAN.

\begin{table}[H]
\centering
\caption{Comparison of errors between training with the same initial sampling and RAR-D adaptive sampling for Example \ref{E5}}\label{table:t15}
\begin{tabular}{c c c c c c}
\toprule
\textbf{Networks} & $\boldsymbol{e_{\Omega_{1}}}$ & $\boldsymbol{e_{\Omega_{2}}}$ & $\boldsymbol{e_{\Gamma}}$  & $\boldsymbol{e_{\partial \Omega_{2}}}$ & $|u-\hat{u}|_{\infty}$ \\
\midrule
PINNs  & $1.409 \times 10^{-3}$ & $3.766 \times 10^{-3}$ & $1.527 \times 10^{-3}$ &  $2.843 \times 10^{-3}$ &  $3.568 \times 10^{-3}$\\

KANs & $1.725 \times 10^{-4}$ & $1.789 \times 10^{-4}$ & $2.599 \times 10^{-4}$ & $2.567 \times 10^{-4}$ &  $5.645 \times 10^{-4}$ \\
PINNs-A & $6.813 \times 10^{-4}$ & $1.689 \times 10^{-3}$ & $8.352 \times 10^{-4}$ &  $1.384 \times 10^{-3}$ &  $1.815 \times 10^{-3}$\\
KANs-A & $4.618 \times 10^{-5}$ & $3.748 \times 10^{-5}$ & $5.695 \times 10^{-5}$ & $6.805 \times 10^{-5}$  &  $1.676\times 10^{-4}$\\

\bottomrule
\end{tabular}
\end{table}

\begin{table}[H]
\centering
\caption{Comparison of errors between training with the same initial sampling and RAR-D adaptive sampling for Example \ref{E6}.}
\label{table:t18}
\resizebox{\linewidth}{!}{%
\begin{tabular}{c c c c c c c}
\toprule
\textbf{Networks} & $\boldsymbol{e_{\Omega_{1}}}$ & $\boldsymbol{e_{\Omega_{2}}}$ & $\boldsymbol{e_{\Gamma}}$ & $\boldsymbol{e_{\partial \Omega_{1}}}$  & $\boldsymbol{e_{\partial \Omega_{2}}}$ & $|u-\hat{u}|_{\infty}$\\
\midrule
PINNs  & $6.628 \times 10^{-4}$ & $2.821 \times 10^{-3}$ & $1.113 \times 10^{-3}$ &  $1.516 \times 10^{-3}$ & $4.559 \times 10^{-3}$ & $9.426 \times 10^{-3}$\\
KANs   & $2.986 \times 10^{-4}$ & $7.623 \times 10^{-4}$ & $4.098 \times 10^{-4}$ & $7.298 \times 10^{-4}$ & $1.552\times 10^{-3}$ & $2.709 \times 10^{-3}$ \\
PINNs-A& $2.131 \times 10^{-4}$ & $1.326 \times 10^{-3}$ & $6.364 \times 10^{-4}$ & $1.429 \times 10^{-4}$ & $3.295 \times 10^{-3}$ & $6.606 \times 10^{-3}$\\
KANs-A & $8.515 \times 10^{-5}$ & $4.902 \times 10^{-4}$ & $1.522 \times 10^{-4}$ & $1.298 \times 10^{-4}$ & $8.365 \times 10^{-4}$ & $1.906 \times 10^{-3}$ \\
\bottomrule
\end{tabular}%
}
\end{table}

%\begin{table}[H]
%\centering
%\caption{Comparison of errors between training with the same initial sampling and RAR-D adaptive sampling for Example \ref{E6}.}\label{table:t18}
%\begin{tabular}{c c c c c c c}
%\toprule
%\textbf{Networks} & $\boldsymbol{e_{\Omega_{1}}}$ & $\boldsymbol{e_{\Omega_{2}}}$ & $\boldsymbol{e_{\Gamma}}$ & $\boldsymbol{e_{\partial \Omega_{1}}}$  & $\boldsymbol{e_{\partial \Omega_{2}}}$ & $|u-\hat{u}|_{\infty}$\\
%\midrule
%PINNs  & $6.628 \times 10^{-4}$ & $2.821 \times 10^{-3}$ & $1.113 \times 10^{-3}$ &  $1.516 \times 10^{-3}$ & $4.559 \times 10^{-3}$ & $9.426 \times 10^{-3}$\\

%KANs & $2.986 \times 10^{-4}$ & $7.623 \times 10^{-4}$ & $4.098 \times 10^{-4}$ & $7.298 \times 10^{-4}$ & $1.552\times 10^{-3}$ & $2.709 \times 10^{-3}$ \\
%PINNs-A & $2.131 \times 10^{-4}$ & $1.326 \times 10^{-3}$ & $6.364 \times 10^{-4}$ & $1.429 \times 10^{-4}$ & $3.295 \times 10^{-3}$ & $6.606 \times 10^{-3}$\\
%KANs-A & $8.515 \times 10^{-5}$ & $4.902 \times 10^{-4}$ & $1.522 \times 10^{-4}$ & $1.298 \times 10^{-4}$ & $8.365 \times 10^{-4}$ & $1.906 \times 10^{-3}$ \\

%\bottomrule
%\end{tabular}
%\end{table}

Selecting $k=2$ and $c=0$. In the adaptive sampling setup,  both networks of PINN and KAN retrained using the Adam for 40000 steps. We generate 10000 initial sampling points due to the larger domain of Example \ref{E5}, and the model is first trained for 20000 steps. Following this, adaptive resampling (RAR-D) is applied after every 2000 iterations.(Example \ref{E5}). Similarly, for Example \ref{E6}, we first generate 5000 initial sampling points and train the model for 20000 steps, subsequently applying adaptive resampling (RAR-D) after every 2000 iterations.

Tables \ref{table:t15} and \ref{table:t18} display the relative $\ell_2$ error of the same initial sampling training and the corresponding adaptive sampling training for Examples \ref{E5} and \ref{E6} respectively.  
As concluded in Example \ref{E4}, we can clearly see that the adaptive algorithm positively improves the performance of both PINN and KAN, while KAN can derive better approximate results than PINN.

\section{Conclusions} 
This paper focused on solving elliptic interface problems using Physics-Informed Neural Networks (PINNs) and Kolmogorov-Arnold Networks (KANs). We adopt a dual-network framework for both PINNs and KANs. Such network admits greater flexibility in addressing the jump conditions at the interface. Through several numerical experiments, we demonstrate that KANs exhibit superior performance in terms of approximate accuracy and error distribution. To further enhance the efficiency of the networks, we incorporate them with the Residual-based Adaptive Refinement with Diversity (RAR-D) sampling strategy. This combination dynamically adjusts the sampling points during the training process, achieving better network performance. The numerical results indicated that KANs not only achieve better approximate solutions but also only require smaller models compared to PINNs, making them a promising alternative for solving elliptic interface problems. This validates that, compared to neural networks based on traditional MLPs, KANs exhibit superior accuracy and interpretability in solving the elliptic interface problems. %Subsequently, the performance of PINNs and KANs in solving higher-dimensional problems can be compared. It also suggests the significant potential of KANs in solving partial differential equations.


\begin{thebibliography}{00}


\bibitem{1}
  Hou T. Y., Li Z., Osher S., and Zhao H.,
  \textit{A hybrid method for moving interface problems with application to the Hele--Shaw flow},
  J. Comput. Phys., 134 (2), pp.~236--252, 1997.

\bibitem{6}
  Lin T. and Wang J.,
  \textit{An immersed finite element electric field solver for ion optics modeling},
  AIAA, 2002.

\bibitem{4}
  Ji H., Zhang Q., Wang Q., and Xie Y.,
  \textit{A partially penalised immersed finite element method for elliptic interface problems with non-homogeneous jump conditions},
  East Asia J. Appl. Math., 8, pp.~1--23, 2018.

%\bibitem{EABE_Rafiezadeh2013}
%Rafiezadeh K., and Ataie-Ashtiani B.,
%\textit{Seepage analysis in multi-domain general anisotropic media by three-dimensional boundary elements},
%Eng. Anal. Bound. Elem., 37, pp.~527--541, 2013.

%\bibitem{EABE_Shiah2021}
%Shiah Y.C., and Hsiao Y.F.,
%\textit{3D analysis of heat conduction in anisotropic composites with thin interstitial layers and imperfect interfaces},
%Eng. Anal. Bound. Elem., 122, pp.~100--114, 2021.

\bibitem{8}
  Ciarlet P. G.,
  \textit{The Finite Element Method for Elliptic Problems},
  North-Holland, 1978.

%\bibitem{11}
%  LeVeque R. J. and Li Z.,
%  \textit{The immersed interface method for elliptic equations with discontinuous coefficients and singular sources},
%  SIAM J. Numer. Anal., 31 (4), pp.~1019--1044, 1994.

\bibitem{Ji_immersed_JSC}
Ji H. and Li Z.,
\textit{An immersed finite element method for anisotropic elliptic interface problems with nonhomogeneous jump conditions}
J. Sci. Comput., 106 (27), 2026.


\bibitem{EABE_Mu2024}
Mu R., Song L., and Qin Q.,
\textit{A meshless method based on the generalized finite difference method for 2D and 3D anisotropic elliptic interface problems},
Eng. Anal. Bound. Elem., 163, pp.~505--516, 2024.

\bibitem{Chu_MG_IFE_JSC}
Chu H., Song Y., Ji H. and Cai Y.,
\textit{Multigrid algorithm for immersed finite element discretizations of elliptic interface problems},
J. Sci. Comput., 98 (26), 2024.

\bibitem{14}
  Raissi M., Perdikaris P., and Karniadakis G. E.,
  \textit{Physics-informed neural networks: A deep learning framework for solving forward and inverse problems involving nonlinear partial differential equations},
  J. Comput. Phys., 378, pp.~686--707, 2019.

\bibitem{Cuomo_PINNs_JSC}
Cuomo S., Schiano Di Cola V., Giampaolo F., Rozza G., Raissi M. and Piccialli F.,
\textit{Scientific machine learning through physics--informed neural networks: where we are and what’s next},
J. Sci. Comput., 92 (88), 2022.

%\bibitem{EABE_ePINN_interface_crack_2024}
%Gu Y., Xie L., Qu W., and Zhao S.,
%\textit{Interface crack analysis in 2D bounded dissimilar materials using an enriched physics-informed neural networks},
%Eng. Anal. Bound. Elem., 163, pp.~465--473, 2024.

\bibitem{18}
  Li W., Xiang X., and Xu Y.,
  \textit{Deep domain decomposition method: Elliptic problems},
  \textit{Math. Sci. Mach. Learn.}, pp.~269--286, 2020.

\bibitem{JSC_DDMPINNs}
Wu W., Feng X. and Xu H.,
\textit{Improved Deep Neural Networks with Domain Decomposition in Solving Partial Differential Equations},
J. Sci. Comput., 93 (20), 2022.


\bibitem{19}
  Jagtap A. D., Kharazmi E., and Karniadakis G. E.,
  \textit{Conservative physics-informed neural networks on discrete domains for conservation laws: Applications to forward and inverse problems},
  Comput. Methods Appl. Mech. Eng., 365, 113028, 2020.

\bibitem{21}
  Jagtap A. D. and Karniadakis G. E.,
  \textit{Extended physics-informed neural networks (XPINNs): A generalized space-time domain decomposition based deep learning framework for nonlinear partial differential equations},
  Commun. Comput. Phys., 28 (5), 2020.

\bibitem{20}
  Sarma A. K., Roy S., Annavarapu C., Roy P., and Jagannathan S.,
  \textit{Interface PINNs (I-PINNs): A physics-informed neural networks framework for interface problems},
  Comput. Methods Appl. Mech. Eng., 429, 117135, 2024.

\bibitem{22}
  Haykin S.,
  \textit{Neural Networks: A Comprehensive Foundation},
  Prentice Hall PTR, 1998.

\bibitem{24}
  Hornik K., Stinchcombe M., and White H.,
  \textit{Multilayer feedforward networks are universal approximators},
  Neural Netw., 2 (5), pp.~359--366, 1989.

\bibitem{Liu_Discontinuity_JSC}
Liu L., Liu S., Xie H., Xiong F., Yu T., Xiao M., Liu L. and Yong H.,
\textit{Discontinuity computing using physics-informed neural networks},
J. Sci. Comput., 98 (22), 2024.

\bibitem{27}
  Liu Z., Wang Y., Vaidya S., R{\"u}hle F., Halverson J., Solja{\v{c}}i{\'c} M., Hou T. Y., and Tegmark M.,
  \textit{KAN: Kolmogorov-Arnold Networks},
  arXiv:2404.19756, 2024.

\bibitem{28}
  Kolmogorov A. N.,
  \textit{On the Representation of Continuous Functions of Several Variables by Superpositions of Continuous Functions of a Smaller Number of Variables},
  American Mathematical Society, 1961.

\bibitem{29}
  Braun J. and Griebel M.,
  \textit{On a constructive proof of Kolmogorov's superposition theorem},
  Constr. Approx., 30, pp.~653--675, 2009.

\bibitem{30}
  Abueidda D. W., Pantidis P., and Mobasher M. E.,
  \textit{Deepokan: Deep operator network based on Kolmogorov Arnold Networks for mechanics problems},
  Comput. Methods Appl. Mech. Eng., 436, 117699, 2025.

\bibitem{31}
  Wang Y., Sun J., Bai J., Anitescu C., Eshaghi M. S., Zhuang X., Rabczuk T., and Liu Y.,
  \textit{Kolmogorov Arnold Informed neural network: A physics-informed deep learning framework for solving PDEs based on Kolmogorov Arnold Networks},
  arXiv:2406.11045, 2024.

\bibitem{33}
  Howard A. A., Jacob B., Murphy S. H., Heinlein A., and Stinis P.,
  \textit{Finite basis Kolmogorov-Arnold networks: domain decomposition for data-driven and physics-informed problems},
  arXiv:2406.19662, 2024.

\bibitem{35}
  So C. C. and Yung S. P.,
  \textit{Higher-order-ReLU-KANs (HRKANs) for solving physics-informed neural networks (PINNs) more accurately, robustly and faster},
  arXiv:2409.14248, 2024.

\bibitem{40}
  Wu C., Zhu M., Tan Q., Kartha Y., and Lu L.,
  \textit{A comprehensive study of non-adaptive and residual-based adaptive sampling for physics-informed neural networks},
  Comput. Methods Appl. Mech. Eng., 403, 115671, 2023.

\bibitem{41}
  Mao Z. and Meng X.,
  \textit{Physics-informed neural networks with residual/gradient-based adaptive sampling methods for solving partial differential equations with sharp solutions},
  Appl. Math. Mech., 44 (7), pp.~1069--1084, 2023.

\bibitem{42}
  Gao Z., Yan L., and Zhou T.,
  \textit{Failure-informed adaptive sampling for PINNs},
  SIAM J. Sci. Comput., 45 (4), pp.~A1971--A1994, 2023.

\bibitem{48}
  Samek W., Montavon G., Lapuschkin S., Anders C. J., and M{\"u}ller K.-R.,
  \textit{Explaining deep neural networks and beyond: A review of methods and applications},
  Proc. IEEE, 109 (3), pp.~247--278, 2021.

\bibitem{45}
  Baydin A. G., Pearlmutter B. A., Radul A. A., and Siskind J. M.,
  \textit{Automatic differentiation in machine learning: a survey},
  J. Mach. Learn. Res., 18 (153), pp.~1--43, 2018.

\end{thebibliography}



\end{document}